\documentclass{article}
\usepackage[utf8]{inputenc}
\usepackage{tabularx}
\usepackage{amsmath}
\usepackage{amssymb}
\usepackage{amsthm}
\usepackage{geometry}
\usepackage{enumitem}
\usepackage[T1]{fontenc}
\usepackage[dvipsnames]{xcolor}
\title{Problems in Modern Galois Theory\footnote{The abstract of this report was given as a lecture at the 1932 International Congress of Mathematics in Zurich.}}
\author{Nikolai Chebotaryov, Translation by Yonathan Stone\thanks{The translator would like to thank Jesse Wolfson for suggesting this and for providing helpful feedback on earlier drafts. In addition, this translation was supported in part by NSF Grant DMS-1944862.}}
\date{First published in 1932, translated in 2020}

\begin{document}
\renewcommand{\thesection}{\S\arabic{section}}
\maketitle

\nocite{*}
Given the recent centennial of Evariste Galois' death, I am presented with the opportunity to present the current state of his most important creation, known by the name of ``Galois theory.''  At the same time, I will try to make a handful of predictions regarding future developments in Galois theory.  The original goal in the development of Galois theory, namely the problem of representing the roots of algebraic equations using radical expressions, was very nearly reached by both Galois and his early followers.  The main tools used by Galois in his investigations, the description of algebraic number fields by their corresponding groups, managed to demonstrate its power additionally for distant branches of mathematical analysis.  In this vein, new topics such as ``Riemann surfaces,'' ``automorphic functions,'' ``continuous transformation groups,'' and so on have been created.

Besides this there are also new problems in classical Galois theory itself, which in turn prompt a deeper analysis of Galois' theory's core ideas.    The problem of finding equations with prescribed [Galois] groups has required studying the theory of generalized rational functions (the problem of Lüroth-Steinitz.)  A generalization of the problem of solving by radicals, namely the problem of Klein forms, has linked the theory of finite groups with the theory of continuous groups.  Connecting both of these branches of group theory is ``Linear group theory,'' a true bridge.

In this report, I would like to present the current machinery available, which may aid the mastery of problems in Galois theory.  In doing so I will broaden the definition of ``Galois theory'' slightly beyond describing the usual application of group theory to algebraic equations. That is, I will include all questions juxtaposing the notions of ``the rational'' and the ``the algebraic irrational.''  Among this are some very nice results concerning algebraic function fields of several variables, which until now have only been facilitated using methods in algebraic geometry.  We can thank the old German and the Italian geometers for solving some of these problems with the help of algebraic geometry, whereas methods for solving such problems are not available to us algebraists.

I have based the selection of this material primarily on my own tastes, and make no claims regarding the objectivity of the following selection:

\S1 Foundations of Galois Theory.\\
\S2 Equations with prescribed [Galois] group.\\
\S3 On the analytic form of prime numbers belonging to a prescribed permutation class.\\
\S4 The resolvent problem.\\
\S5 Further questions in the general theory of fields.
\section{ Foundations of Galois Theory}
\numberwithin{equation}{section}
\renewcommand{\thesection}{\arabic{section}}
\begin{enumerate}
    \item One can separate the work which deals with the foundations of Galois theory into two types.  The first of these includes the tasks which seek new paths to justify the classical Galois theory, while those of the second kind deepen the definition of the Galois group, facilitating its application to far more areas than is possible using just the classical theory.
    \item Belonging to the first kind, the work of F. Mertens \cite{ref50}, S. Schatunowski \cite{ref62}, and A. Loewy \cite{ref46, ref47} is of particular note.  Mertens develops the definition of the Galois group and proves the associated fundamental theorems without using the notion of the normal field [extension], the Galois resolvent, etc. He builds upon the definition of irreducibility in extended domains. Given an equation $f(x) = 0$ with $x_1$ as a root, he finds a factor $f_1(x;x_1)$ of the polynomial $\frac{f(x)}{x-x_1}$ that is irreducible over $K[x_1]$.  He proceeds to find a $K[x_1,x_2]$-irreducible factor of the polynomial $\frac{f_1(x;x_1)}{x-x_2}$, where $x_2$ is a root of $f_1(x;x_1)$. Continuing this process, he arrives at the family
    \begin{equation}\label{1.1}
      Z_0 = f(x), Z_1 = f_1(x;x_1), Z_2 = f_2(x;x_1,x_2), ..., Z_{n-1} = f_{n-1}(x;x_1,...,x_{n-1})  
    \end{equation}
    of polynomials, which forms a family of \textit{fundamental relations}. A polynomial $\varphi(x_1,x_2,...,x_n)$ is then equal to zero (in the variables $x_i$) if and only if it can be represented in the form
    \begin{equation}\label{1.2}
        P_0\cdot Z_0(x_1) + P_1\cdot Z_1(x_2) + ... + P_{n-1}\cdot Z_{n-1}(x_n),
    \end{equation}
    where the $P_i$'s are polynomials.  The Galois group is thus defined to consist of those permutations that transform each of the $Z_i(x_{i+1})$ into a polynomial of the form (\ref{1.2}).  The order of this Galois group is equal to the product of degrees of the polynomials in (\ref{1.1}).
    \item Schatunowski establishes a theory formally covering the same ideas as Mertens, yet at the same time assumes a much more general point of view.  His work builds upon Kronecker's idea of functional relations, which, given a root $x_i$ of the equation
    \begin{equation}\label{1.3}
        f(x) = x^n + a_1x^{n-1} + ... + a_{n-1}x + a_n = 0,
    \end{equation}
    consists of considering any quantity depending on $x_i$ as a function of an undetermined variable, while taking congruence modulo $f(x)$ as the equals sign.  Since one considers functions of several roots of equation (\ref{1.3}) in Galois theory, Schatunowski aims to construct a family of functional relations in the variables $x_1,...,x_n$, such that the quotient according to these relations is isomorphic to the algebraic number field generated by the roots of equation (\ref{1.3}).  The ideal generated by
    \begin{align}\label{1.4}
        \psi_1 = x_1 + x_2 + ... + x_n + a_1, &\psi_2 = x_1x_2 + ... + x_{n-1}x_n - a_2,..., \\ 
        &\psi_n = x_1x_2...x_n - (-1)^na_n \nonumber
    \end{align}
    is not suited for this, since the $f(x_i)$ are not $\equiv 0$ modulo this ideal.  The ideal generated by
    \begin{equation}\label{1.5}
        f(x_1),f(x_2),...,f(x_n)
    \end{equation}
    also fails to serve this purpose, since the statement
    \[\psi_i \equiv 0 (\text{ mod } f(x_1),f(x_2),...,f(x_n))\]
    fails to hold, whereas $\psi_i V \equiv 0 (\text{ mod } f(x_1),f(x_2),...,f(x_n))$ does hold, where $V$ denotes the Vandermonde determinant of $x_1,...,x_n$.  In accordance with modern convention we say that (\ref{1.4}) and (\ref{1.5}) are not prime ideals. \\
    Schatunowski constructs the desired system of relations by starting with what he calls Cauchy relations.  One obtains these relations in the following manner:  As the first relation one takes $f(x_1)$; as the second relation the quotient of the division of $f(x_1)$ by $x_1 - x_2$; as the third relation the quotient of the second relation, ordered by powers of $x_2$, by $x_2 - x_3$, etc. The Cauchy relations then give rise to residue classes which in turn are isomorphic to the corresponding number field if and only if equation (\ref{1.3}) is without affect. [German: \textit{affektlos}]  In this case each of these ideals is irreducible over the ideal preceding it.  If this is not the case Schatunowski considers the general Mertensian relations.  The degrees of these relations provide information regarding the transivity and primitivity relations of equation (\ref{1.3}).  For instance, the [Galois] group is $k$-times transitive whenever the first $k$ relations coincide with those of Cauchy.
    \\
    Schatonowski's work contains the introductory foundations of what is known as the \textit{theory of polynomial ideals} nowadays.  It deals with the rings whose relations are reducible.  Such a ring is called a \textit{``semifield''} if one can adjoin a (finite) number of new relations such that the ring becomes a field.\\
    Paritcularly notable among terms introduced by Schatunowski is the notion of an \textit{extension of the second kind}.  He uses this to refer to rings that are obtained by adjoining a finite number of moduli to a base ring.  In other words, one sets certain nonvanishing elements of the base ring equal to zero.  It may be more natural to refer to these extensions as \textit{quotient rings}, in accordance with the term \textit{quotient group.}  It is known that a field does not admit any extensions of the second kind.  In this, it is assumed that all fields can only have characteristic zero.  However, if we consider the prime characteristic $p$, this fails upon adding the new modulus $p$, and a field remains a field.
    \\
    Assuming the Mertens-Schatunowskian definition of the Galois group, one can easily prove the following theorem, originally due to I. Schur (\cite{ref70}, see also \cite{ref79}):\\
    The [Galois] group of a factor field is a subgroup of the group of the original field.
    \item The Loewyian justification of Galois theory has some overlap with the Mertens-Schatunowskian theory, in particular the avoiding the assumption of normal fields.  Rather than taking a single primitive element, Loewy takes several algebraic elements as possible generators for a field $P$, which he refers to as \textit{conductors}\footnote{Translated from the German ``Dirigenten''.-\textit{YS}}.  The first conductor $\rho_1$ is a root of an equation with coefficients in the domain of rationality; the second conductor is a root of an equation of the type $f(\rho_1;z) = 0$; the third of the type $f(\rho_1,\rho_2;z) = 0$, etc.  Provided all these equations are irreducible, Loewy uses ``transformations of the field $P$'' to refer to the replacement of a conductor $\rho_i$ with a conjugate root belonging to all equations in which $\rho_i$ occurs (which also necessitates the replacement of all subsequent conductors).  He proves that such a transformation disturbs none of the relations between the conductors.  However, it is in fact possible that such a transformation maps some elements of the field $P$ outside of $P$.  The transformations that map elements $P$ to $P$ form a group, which is called the group of \textit{automorphic transformations}. From a group theoretic standpoint this means the following:  A non-normal field $P$ is not at all determined by its [Galois] group; it is determined by its Galois group $\mathfrak{G}$ and the subgroup $\mathfrak{H}$ corresponding to a primitive element of $P$.  The automorphic group of $P$ is isomorphic to $\mathfrak{K}/\mathfrak{H}$, where $\mathfrak{K}$ denotes the normalizer of $\mathfrak{H}$ in $\mathfrak{G}$.  The collection of all transformations of $P$, using the above definition, do not form a group in the familiar sense.  Loewy calls these algebraic structures \textit{mixed groups} and investigates their structural properties \cite{ref47}.  Each mixed group $\mathfrak{T}$ contains a kernel $\mathfrak{G}$ that is defined to be the largest ordinary group contained inside $\mathfrak{T}$ and consists of some of the cosets  (residue classes) of $\mathfrak{G}$. It is important, that a factor group of $\mathfrak{T}$ by \textit{any} (not neccesarily normal) subgroup $\mathfrak{H}$ of $\mathfrak{G}$ remains a mixed group, whose kernel is given by $\mathfrak{N}/\mathfrak{G}$, where $\mathfrak{N}$ denotes the normalizer of $\mathfrak{H}$ in $\mathfrak{G}$.  These facts permit viewing a mixed group as a more adequate picture of a field.
    \\ A similar group construction was introduced by H. Brandt \cite{ref7}, which is called the \textit{Brandt groupoid}.
    \item Before we transition to the work of the second kind, we must set out the modern notion of the term ``Galois group'', which deviates from older notions.  It is difficult for me to say who this new notion originates from.  The older Galois theory views the elements of the Galois group, the \textit{transformations} (or \textit{permutations}), as interchanges among the roots of a generating equation (which may well be reducible), which do not disturb any relations between the roots. On the other hand, the modern Galois theory considers assignments that are simultaneously undergone by all elements of a field $K$, without disturbing the established relations between them.  In other words, each element of a Galois group is a map of a (normal) field $K$ to itself, or, as one says in group theory, an automorphism.  That is, an assignment  of all elements of a field to elements of the same field, such that sums and products and taken to sums and products respectively.
    \item Of utmost importance in Galois theory is the mutual correspondance between subfields of $K$ and subgroups of its Galois group.  More specifically one can formulate this in the following manner (\cite{ref43}; \cite{ref74}, appendix; \cite{ref3}):  One assigns to each subfield $U$ of $K$ the largest subgroup $\mathfrak{H}(U)$ of $\mathfrak{G}$ which fixes all elements of $U$.  On the other hand one can assign to each subgroup $\mathfrak{H}$ of $\mathfrak{G}$ the largest subfield $U(\mathfrak{H})$ of $K$, whose elements remain unchanged by $\mathfrak{H}$.  If $U_1 > U_2$, then $\mathfrak{H}(U_1) < \mathfrak{H}(U_2)$, and vice-versa.  Moreover it holds that
    \begin{equation}\label{1.6}
        U[\mathfrak{H}(U)] > U, \mathfrak{H}[U(\mathfrak{H})] > \mathfrak{H}.
    \end{equation}
    However, one can only readily develop Galois theory, whenever
    \begin{equation}\label{1.7}
        U[\mathfrak{H}(U)] = U,
    \end{equation}
    \begin{equation}\label{1.8}
        \mathfrak{H}[U(\mathfrak{H})] = \mathfrak{H}.
    \end{equation}
    So that (\ref{1.8}) holds, $K$ is required to be finite over its domain of rationality (Krull, \cite{ref43}).\\
    So that (\ref{1.7}) holds, $K$ must be of the first kind over its domain of rationality (Baer-Hasse, \cite{ref74}, appendix).\\
    If $K$ possesses the field of rational numbers as a subfield  (one says that $K$ has characteristic zero), then $K$ is of the first kind in any case.  On the other hand, if $K$ has characteristic $p$ (that is, there exists a prime number $p$, that is equal to zero inside $K$), then $K$ is of the first kind if and only if a generating element of $K$ satisfies an irreducible equation with entirely distinct roots\footnote{This coincides precisely with the modern notion of separability.-\textit{YS}}. \\
    K is finite over its domain of rationality, if there exists a finite number of basis elements, so that each element of $K$ has a linear representation with in terms of the basis elements and coefficients from the domain of rationality.
    \item For the case of $K$ being infinite, Krull \cite{ref43} broadened the fundamental theorem of Galois theory by considering only the \textit{closed} subgroups of $\mathfrak{H}$ in lieu of all subgroups.  By this he means the following.  Given $\gamma$ an element of $\mathfrak{G}$ and $x$ a finite normal subfield of $K$, a mapping of $x$ to itself is induced by $\gamma$.  If both $\gamma$ and $\gamma^*$ induce the same mapping of $x$, we say that $\gamma^*$ is contained in an $x$-\textit{neighborhood} of $\gamma$.  This definition of neighborhood satisfies all the Hausdorff neighborhood axioms:
    \begin{enumerate}[label=\alph*)]
        \item Each element $\gamma$ is contained in each of its own neighborhoods.
        \item The intersection of two neighborhoods of $\gamma$ contains a new neighborhood of $\gamma$.  Moreover, it itself is a neighborhood of $\gamma$, since the intersection of the $x_1-$ and $x_2$-neighborhoods is the $x_3$-neighborhood, where $x_3$ is both finite and normal as it is the union\footnote{The author uses the German word ``\textit{Vereinigung}'', or union, but in this case the notion of a \textit{compositum of fields} seems more appropriate.-\textit{YS}} of $x_1$ and $x_2$.
        \item Given $\delta$ an element of a $x$-neighborhood of $\gamma$, there exists a neighborhood of $\delta$ which is entirely contained inside the $x$-neighborhood of $\gamma$.  Moreover, the $x$-neighborhoods of $\gamma$ and $\delta$ coincide, since they contain the largest number of elements of $\mathfrak{G}$ that act on the field $x$ in precisely the same way.
    \end{enumerate}
    If $\gamma$ and $\delta$ are distinct elements of $\mathfrak{G}$, there exists a neighborhood of $\gamma$ not containing the element $\delta$.  This is because $\gamma$ and $\delta$ being distinct implies the existence of elements in $K$ that behave differently under $\gamma$ and $\delta$.  As each of these elements generates a finite field [extension], it follows that each of the fields corresponds to a different neighborhood of $\gamma$ and $\delta$.
    \\
    This definition permits the definition of \textit{accumulation elements}.  Given $\mathfrak{H}$ a subgroup of $\mathfrak{G}$, an accumulation element of $\mathfrak{H}$ should contain elements of $\mathfrak{H}$ in each of its neighborhoods. An accumulation element of $\mathfrak{H}$ is not necessarily contained in $\mathfrak{H}$ itself. However, in the case where a subgroup $\mathfrak{H}$ of $\mathfrak{G}$ contains all of its accumulation elements, it is called \textit{closed}.  Each group which corresponds to a subfield of $K$ is closed.  Conversely, to each closed group $\mathfrak{H}$ there is a corresponding subfield $U$ of $K$, so that $\mathfrak{H}[U(\mathfrak{H})] = \mathfrak{H}$ holds true.  Thus, in order to be able to broaden the fundamental theorems of Galois theory to infinite fields [extensions], one must only consider those subgroups of $\mathfrak{H}$ which are closed.
    \item The conditions for the presence of the relation (\ref{1.7}) were investigated thoroughly by R. Baer \cite{ref3}.  He found that, in any case, there exists an intermediate field [extension] $S$, which he called the \textit{rigid field} between $K$ and the domain of rationality.  The rigid field is characterized by $K$ being \textit{orderly} (that is (\ref{1.8}) always holds) when one takes $S$ to be the domain of rationality, while all elements of $S$ remain invariant under automorphisms of $K$.  I cannot elaborate on the additional interesting statements of this work here.
    \item It is very difficult to define the Galois group in the case where the field in question $K$ has a higher degree of transcendence than its domain of rationality.  The reason for this lies in the fact that the \textit{universal norm} of such a field (that is the field, which contains all fields conjugate to subfields of $K$) is an infinite field [extension] whose definition is difficult to summarize analytically.  To theoretically construct a group which possesses the desired main properties of the Galois group, one can sketch the following scheme.  Let $x_1,x_2,...,x_n$ be the generating elements of a field $K$, among which certain algebraic relations may be established, which we will denote by (I).  One can determine each subfield $U$ of $K$ in an analogous sense through generating elements $\xi_1,\xi_2,...,\xi_m$, where the $\xi_i$ are expressed rationally through the $x_i$:
    \[\xi_i = x_i(x_1,x_2,...,x_n) \quad \quad (i = 1,2,...,m).\]
    The equations
    \begin{equation}
        \xi_i(x_1,x_2,...,x_n) = \xi_i(y_1,y_2,...,y_n) \quad \quad (i = 1,2,...,m) \tag{II}
    \end{equation}
    determined a new field, whose generators $[x_1,x_2,...,x_n;y_1,y_2,...,y_n;y_1',y_2',...,y_n';...]$ are related by the relations (I) and (II).  This fields can be called the \textit{relative norm} of $K$ (with respect to $U$).  The assignment of $(x_1,...,x_n)$ to $(y_1,...,y_n)$ is designated a transmutation of $K$ or a permutation of its relative norm.  If $U$ runs through all subfields of $K$, then the sought after universal norm is generated by the relative norms.  The compositum of all previously constructed permutations is a group possessing all main features of the Galois group.
    \item One can construct the Galois group of a field of algebraic functions in a slightly different manner.  Instead of considering the functions of the field, one envisions the totality of their coordinate systems, which form a so-called \textit{absolute Riemann surface} (see \cite{ref86}).  In this case each transformation of the Galois group assigns to each value of a function of $K$ a new value, such that the relations between values of different functions are preserved.  As each function is determined by the totality of its values, such transformations also determine the functions to which the given functions are assigned.  The various \textit{monodromy groups} which do not disturb certain domains of rationality are contained in this group.  It can very well occur that a transformation leads several functions of $K$ outside of $K$.  In the case of an independent variable this is reflected when a function determined up to a multiplicative constant by its zeroes and poles [of the form]
    \[f \simeq \frac{\mathfrak{p}_1'\mathfrak{p}_2'...\mathfrak{p}_m'}{\mathfrak{p}_1\mathfrak{p}_2...\mathfrak{p}_m}\]
    is assigned to a product
    \[\frac{\overline{\mathfrak{p}_1'}\;\overline{\mathfrak{p}_2'}\;...\overline{\mathfrak{p}_m'}}{\overline{\mathfrak{p}_1}\;\overline{\mathfrak{p}_2}\;...\overline{\mathfrak{p}_m}}\]
    where the numerator and denominator are contained in different ideal classes.  Each transformation assigning divisors to equivalent divisors  belongs to the so-called \textit{group of intrinsic transformations} \cite{ref38}, which play an analogous role to the group of \textit{automorphic transformations} introduced by A. Loewy \cite{ref47}.
    \item Of importance is the theory of algebraic functions of a group, which is closely related to the Galois group which has just been defined.
    \\
    Let $\mathfrak{u}_1^{\mathfrak{p},\mathfrak{p}'},\mathfrak{u}_2^{\mathfrak{p},\mathfrak{p}'},...,\mathfrak{u}_p^{\mathfrak{p},\mathfrak{p}'}$ be the linearly independent Abelian integrals of the first kind defined on the Riemann surface of $K$.  The Jacobi inversion problem consists in the solution of the system of equations
    \begin{equation}\label{1.9}
        \mathfrak{u}_i^{\mathfrak{p}_1,\mathfrak{p}_1'} + \mathfrak{u}_i^{\mathfrak{p}_2,\mathfrak{p}_2'} + ... + \mathfrak{u}_i^{\mathfrak{p}_p,\mathfrak{p}_p'} \equiv \mathfrak{v}_i \quad (i = 1,2,...,p),
    \end{equation}
    where the lower limits $\mathfrak{p}_i$ are given and the upper limits $\mathfrak{p}_i'$ are sought for and the congruences are taking from the period systems as relations.  When none of the points $\mathfrak{p}_i$ are collinear, this problem is clearly solvable (\cite{ref54};\cite{ref4}).  If we take the $\mathfrak{p}_i$, $\mathfrak{p}_i'$ to be coordinates of the points $P$, $P'$ in a $p$-dimensional space, then each coordinate system of the parameters $\mathfrak{v}_i$ determines a transformation assigning each point $P$ to a certain point $P'$ (in this case one must not regard points differing by the order of their coordinate values as different, so that these points may be defined unambiguously using the values of the symmetric functions of $\mathfrak{p}_i$).  The group of these transformations is \textit{locally isomorphic}\footnote{The author uses the term \textit{im Kleinen isomorph}, where the adjective \textit{im Kleinen} seems to be a weaker statement that local when used to qualify other definitions.  Whether this is the same as locally isomorphic seems uncertain and as such any clarification would be greatly appreciated.-\textit{YS}} to the $p$-parameter group of translations and admits an analytical expression in terms of the addition formulas of the Abelian functions.  It is a subgroup of an extended Galois group which transforms the coordinates of the point $P$ independently of each other.  We will refer to this group as the \textit{Jacobian group} in the sequel.
\end{enumerate}
\renewcommand{\thesection}{\S\arabic{section}}
\section{Equations with prescribed [Galois] group}
\renewcommand{\thesection}{\arabic{section}}
\begin{enumerate}
    \item The problem of finding equations with prescribed [Galois] groups belongs to the most important problems of modern Galois theory and to this day has not been solved.  It can be summarized in the following three ways:
    \begin{enumerate}[label = \Roman*.]
        \item One finds any equations, whose [Galois] groups are isomorphic to a given group $\mathfrak{G}$.
        \item One finds the most general parametric form of the coefficients of an equation whose [Galois] group is isomorphic to $\mathfrak{G}$ or a subgroup of $\mathfrak{G}$.  Being able to present the coefficients in this form should provide a necessary and sufficient condition for the [Galois] group of the equation to be isomorphic with $\mathfrak{G}$ or one of its subgroups.
        \item One outlines a procedure for the determination of equations, whose [Galois] group is isomorphic to $\mathfrak{G}$.  This procedure should produce all equations of this kind if continued sufficiently.
    \end{enumerate}
    \item Problem II always admits a solution whenever the generalized L{\"u}roth theorem (also called the \textit{rational minimal basis} theorem) holds for a given group $\mathfrak{G}$.  This theorem may be formulated in the following manner:
    \\
    For $K_n(x_1,x_2,...,x_n)$ the field of rational functions of the variables $x_1,x_2,...,x_n$, every subfield of $K_n(x_1,x_2,...,x_n)$ is isomorphic to $K_m(x_1,x_2,...,x_m)$ $(m \leq n)$.  (One can also say: this field is \textit{purely transcendental})\\
    This theorem was proven in the case $n = 1$ by L{\"u}roth \cite{ref48}.  The proof for $n = 2$ was found by Castelnuovo \cite{ref15}.  Both G. Fano \cite{ref21} and F. Enriques \cite{ref20} found a counterexample to the case $n = 3$.  However, the full scope of this theorem is not required for the solving Problem II.  Restricting ourselves to the case where the subfield in question contains the field of elementary symmetric functions in $x_1,x_2,...,x_n$, the validity of the theorem becomes fully equivalent to the solvability of task II.  Until now the validity of this theorem ``in a narrower setting'' remains open (see \cite{ref74}, remark of B.L. Van der Waerden).  If this also does not hold true in general, one can determine whether any abstractly given finite group is ``L{\"u}rothian'' or not.  That is, whether the field $K(a_1,a_2,...,a_n;\varphi)$ is purely transcendental, where $\mathfrak{G}$ is represented as the permutation group of $x_1,x_2,...,x_n$, $a_1,a_2,...,a_n$ denote the elementary symmetric functions of $x_1,x_2,...,x_n$, and $\varphi$ is a function of $x_1,x_2,...,x_n$ belonging to $\mathfrak{G}$.
    \item E. Noether \cite{ref55} deduced Problem I from Problem II by using Hilbert's irreducibility theorem, which states for any irreducible polynomial one can choose numerical values for some of the variables so that the polynomial is irreducible in the remaining variables.
    \\
    One can tighten this result by not just deducing Problem I, but also Problem III from Problem II.  For that one uses the following procedure specified by M. Bauer \cite{ref5},\cite{ref6}.  It is known that the Galois group of an algebraic equation
    \begin{equation}\label{2.1}
        x^n + a_1x^{n-1} + ... + a_{n-1}x + a_n = 0
    \end{equation}
    is transformed into one of its subgroups when one regards the values of its domain of rationality not absolutely, but modulo a prime number (or a prime ideal).  That is, one replaces the domain of rationality with one of its factor rings \cite{ref16}\cite{ref70}\cite{ref79}.  Then again, it also known that the Galois group of a field of characteristic $p$ [German: \textit{Primzahlmodulkongruenz}, or prime number modular congruence] is cyclic, and that a generating permutation of the latter group consists of cycles, whose orders are equal to the irreducible components of our congruence.  As a result, if
    \begin{equation}\label{2.2}
    f(x) \equiv X_{n_1}^{(p)}X_{n_2}^{(p)}...X_{n_k}^{(p)} \quad (\text{ mod } p),
    \end{equation}
    holds, where $X_{n_i}^{(p)}$ denotes a modulo $p$ irreducible polynomial of degree $n_i$ ($n_1 + n_2 + ... + n_k = n$), then the [Galois] group of (\ref{2.1}) contains a permutation whose cycles have order $n_1,n_2,...,n_k$.
    \item We assume that the field $K(a_1,a_2,...,a_n;\varphi)$ is purely transcendental (see the definitions in no. 2). There thus exist rational functions $\pi_1,\pi_2,...,\pi_n$ of $a_1,a_2,...,a_n;\varphi$ so that $a_1,a_2,...,a_n;\varphi$ can themselves be expressed in terms of the $\pi_i$.  Next let $\mathfrak{G}_1,\mathfrak{G}_2,...,\mathfrak{G}_s$ be a system of subgroups of $\mathfrak{G}$, so that each proper subgroup of $\mathfrak{G}$ is a subgroup of at least one of the $\mathfrak{G}_i$.  Such a system can be certainly be constructed, for instance by taking all proper subgroups of $\mathfrak{G}$ as the $\mathfrak{G}_i$.\\
    Let $\varphi_i$ be a function of $x_1,x_2,...,x_n$ belonging to $\mathfrak{G}_i$ $(i = 1,2,..,s)$, and let $F_i(x_i)$ be the polynomial of smallest degree, whose coefficients are rational functions of $\pi_1,\pi_2,...,\pi_n$ and for which $\varphi_i$ is a root ($i = 1,2,...,s$).  One can see that the degree of $F_i(z_i)$ is equal to the index $(\mathfrak{G}:\mathfrak{G}_i)$ ($i = 1,2,...,s)$.\\
    The [Galois] group of the equation $F_i(z_i) = 0$, as a transitive permutation group, contains a permutation $\overline{S_i}$, which fixes no indices.  This permutation corresponds to at least one permutation in $\mathfrak{G}$.  Let $S_i$ be such a permutation, and let $n_1,n_2,...,n_k$ be the orders of its cycles.  Take any prime number $p_i > n-2$ and let
    \begin{equation}\label{2.3}
        f(x) \equiv X_{n_1}^{(p_i)}X_{n_2}^{(p_i)}...X_{n_k}^{(p_i)} \quad (\text{ mod } p_i),
    \end{equation}
    where $X_{n_j}^{(p_i)}$ denotes a modulo $p_i$ irreducible polynomial of degree $n_j$.  This determines the congruence classes modulo $p_i$ for the $a_i$.  Using the expressions for the coefficients of the equation $F(z) = 0$, and we obtain the congruence
    \[F(z) \equiv 0 \quad (\text{ mod } p_i),\]
    which surely has one or more rational roots.  Let $\varphi_i$ be one of these roots.  Then $\varphi_i$ is transformed into the other rational roots by means of some permutations $\Sigma_1,\Sigma_2,...,\Sigma_\nu$ of the symmetric permutation group $\mathfrak{S}$ of $x_1,...,x_n$.\\
    By fixing (\ref{2.3}), the ``broad'' class of permutations containing $p_i$ is determined, that is the totality of all permutations similar to $S_i$ within $\mathfrak{S}$.  However, by fixing the congruence class of $\varphi_i$, we determine a \textit{division} \footnote{From the German \textit{Abteilung}.  According to \textit{Pioneers of Representation Theory: Frobenius, Burnside, Schur, and Brauer}
by Charles W. Curtis this is a term a term introduced by Frobenius.  A \textit{division} of a [Galois] group is the union of conjugacy classes of a given group element and all of the powers of this element coprime to its order.-\textit{YS}}.  If $\mathfrak{U}$ is a division of $S_i$, then
    \[\Sigma_1^{-1}\mathfrak{U}\Sigma_1^{-1}, \Sigma_2^{-1}\mathfrak{U}\Sigma_2^{-1},..., \Sigma_\nu^{-1}\mathfrak{U}\Sigma_\nu^{-1}\]
    are precisely those divisions of $\mathfrak{G}$ possessing the cycle type of $S_i$.  One of these division must correspond to the division of a permutation $\overline{S}_i$ of our choosing.  Plugging the values of $\varphi_i^{\Sigma_j} (j = 1,2,...,\nu)$ into the expressions for the coefficients of $F_i(z_i) \equiv 0$, then at least one of the resulting congruences $F_i(z_i) \equiv 0 (\text{ mod } p_i)$ corresponds to $\overline{S}_i$ and thus has no rational roots.  In other words $p_i$ belongs to a permutation class of $K(x_1,x_2,...,x_n)$ which is not contained in $\mathfrak{G}_i$.\\
    Taking $i = 1,2,...,s$, we obtain the congruence classes modulo $P = p_1p_2...p_s$ for $a_1,a_2,...,a_n;\varphi$ and thus for $\pi_1,\pi_2,...,\pi_n$.  Fixing the $\pi_i$ within the congruence classes we have just determined and plugging them into the equation $f(x) = 0$, we obtain precisely $\mathfrak{G}$ as the group of this equation.  On the one hand, it is contained in $\mathfrak{G}$ by virtue of the parametric expressions of the $a_i$, but on the hand is not contained in any of the groups $\mathfrak{G}_1, \mathfrak{G}_2,...,\mathfrak{G}_s$ due to the established congruence conditions.
    \item In particular, if one wants to construct equations without affect, one can follow a process due to M. Bauer.  Take three arbitrary prime numbers $p,q,r$ ($r \geq n -2$) and impose the following three congruence conditions on the polynomial $f(x)$:
    \begin{align*}
        f(x) &\equiv X_n^{(p)} \quad \text{(mod) $p$},\\
        f(x) &\equiv X_{n-1}^{(q)}(x-b) \quad \text{(mod) $q$}, \\
        f(x) &\equiv X_2^{(r)}(x-b_1)(x-b_2)...(x-b_{n-2})\quad \text{(mod) $r$}.
    \end{align*}
    The [Galois] group of the equation $f(x) = 0$ thus contains an $n$-cycle, an $(n-1)$-cycle and a transposition and thus the symmetric group \cite{ref5}\cite{ref6}\cite{ref79}.
    \item The procedure outlined in No. 4 allows the construction of all possible equations with [Galois] group $\mathfrak{G}$ provided one continues the process sufficiently.  This follows from the following result of Frobenius \cite{ref23}:\\
    If the [Galois] group of the equation $f(x) = 0$ contains a permutation with the cycles of length $n_1,n_2,...,n_k$ ($\Sigma n_i = n)$, then there exist infinitely many prime numbers $p$ such that the congruence $f(x) \equiv 0 \text{ (mod $p$)}$ splits into irreducible polynomials of degree $n_1,n_2,...,n_k$. \\
    The somewhat vague term ``all equations'' can be made precise by giving ourselves the problem of finding all equations with [Galois] group $\mathfrak{G}$, such that their coefficients don't exceed a certain bound.  To this end one must tighten the result of Frobenius in the following sense:
    \\
    One finds bounds under which a prescribed number of prime numbers with the desired properties are located.\\
    Such bounds have been provided by L. Kronecker \cite{ref42} and F. Mertens \cite{ref49} in the case of arithmetic progressions.  I have carried out this estimate for the problem of Frobenius \cite{ref81}.  The result is as follows. If
    \begin{equation} \label{2.4}
    x = \text{Max}\left\{2\left(\frac{2A_d}{g_dh_d} + 2W\right)^{\tfrac{1}{a}}\right\}
    \end{equation}
    for all $d \mid f$,  then the interval $(1,x)$ surely contains $V$ prime numbers belonging to the division of $S$.  The constants $A_d, g_d, h_d, W$ depend on certain subfields of $K$ and on the number $V$.  In order to express this bound explicitly in terms of the coefficients of equation (\ref{2.1}), it is neccesary to estimate certain constants of K.  R. Remak \cite{ref61} very recently constructed an upper and a lower bound for the regulator of a field, which is particularly important for estimating formula (\ref{2.4}).
    \item If we seek only a solution to Problem I the solvability of L{\"u}roth's problem is unneccesary.  If $\varphi$ is a function belonging to $\mathfrak{G}$ which satisfies the equation $F(\varphi) = 0$, then the coefficients of the polynomial $F(\varphi)$ are given by rational functions of the coefficients $a_1,a_2,...,a_n$ of the polynomial $f(x).$  Taking $f(x) \equiv X_{n_1}^{(p_j)}X_{n_2}^{(p_j)}...X_{n_k}^{(p_j)} \quad \text{(mod $p_j$)}$, we can determine the $a_i$ modulo $p_j$.  Plugging these values into the congruence $F(z) \equiv 0 \quad \text{(mod $p_j$)}$, we obtain at least one rational root.  If the congruence has multiple roots, we choose the one that corresponds to the desired division.  Letting $j$ run through $1,2,...,s$, we can determine $a_1,a_2,...,a_n;\varphi$ moduli $p_1,p_2,...,p_s$ and thus modulo $P = p_1p_2...p_s$, and can thus be represented in the form $a_i = a_i^{(0)} + Pt_i$, $\varphi = \varpi_0 + Pu$, where $a_1^{(0)},a_2^{(0)},...,a_n^{(0)};\varphi_0$ denote constants.  Substituting this into the equation $F(\varphi) = 0$, we obtain the Diophantine equation $\Phi(t_1,t_2,...,t_n;u) = 0$.  The Problem now follows from the solution of this equation.  Note that this equation possesses the following properties:
    \begin{enumerate}[label = \roman*)]
        \item It can always be solved using fractions. That is, one can replace the $a_i$ with the elementary symmetric functions of the $n$ arbitrary rational numbers.  This means the solution has $n$ ``degrees of freedom.''
        \item It can always be solved using $p$-adic numbers, where the prime number $p$ can be chosen arbitrarily.
    \end{enumerate}
    The solvability of the equation $\Phi = 0$ does not depend on the choice of function belonging to $\mathfrak{G}$.
    \item Much work has been devoted to the solution of the problems in question in the cases of a few special groups.  In the first place we have the work of D. Hilbert \cite{ref35}, in which Problem I is solved for symmetric and alternating groups of any order using the irreducibility theorem.  The genuine  construction of equations with alternating groups was recently carried out to near completion by I. Schur \cite{ref72}.  That is, given $n \equiv 0 \text{mod $4$}$, Schur shows that the equation
    \[E_n(x) = 1 + \frac{x}{1!} + \frac{x^2}{2!} + ... + \frac{x^n}{n!} = 0\]
    possesses the alternating group as its Galois group.  On the other hand he shows (Crelle: 165, 1931) that when $n \equiv 1 \text{(mod) 2}$, the equation
    \[1 - {n\choose1}\frac{x}{2!} + {n \choose 2}\frac{x^3}{3!} - ... + (-1)^n \frac{x^n}{(n+1)!} = 0\]
    also has the alternating group as its Galois group.
    \item In the solution of Problem II it is important to solve L{\"u}roth's problem using  rational functions with \textit{rational coefficients}.  This question in the case of solvable groups of prime degree is addressed in the work of S. Breuer \cite{ref10}\cite{ref11} and Ph. Furtw{\"a}ngler \cite{ref27}.  Furtw{\"a}ngler proposed the following sufficient condition in the case of cyclic groups of prime degree $p$:
    
    The problem permits a solution whenever it is possible to construct an integer number system $e_0, e_1,...,e_{p-2}$, such that $1)$ the Hankel determinant is
    \[\left| \begin{matrix} e_0, & e_1, & ..., & e_{p-2} \\
    e_{p-2}, & e_0, & ..., & e_{p-3} \\
    ... & ... & ... & ... \\
    e_1, & e_2, & ..., & e_0
    \end{matrix}\right| = \pm p,\]
    and $2)$ the congruences $\sum\limits_{i=0}^{p-2} e_ig^i \equiv 0 \text{(mod $p$)}$ hold, where $g$ is a primitive root of $p$.
    \\
    This criterion is not always satisfied, for instance in the case of $p = 47$. Later, Furtw{\"a}ngler gives some general instructions for the construction of rational minimal bases of metacyclic groups.\\
    Breuer derived several similar criteria by using his theorem about the decomposition of a field of rational functions of $n$ variables into two subfields, where one depends on the full metacyclic functions of the generating variables.
    \item Problems I and III have been extended to the setting of relative fields.  First of all, the solution of problem III can be regarded as solved for relative fields when a known rational basis contains certain roots of unity in the coefficients rather than being rational-valued (see for instance E. Fischer, \cite{ref22}).  This is however not an extension of the problem.  A solution to problems I and III can only be considered a satisfactory extension to relative fields if we are capable of also saying something about the absolute Galois group.\\
    One can call the work of O. Ore \cite{ref57} and H. Hasse \cite{ref30}\cite{ref31} as the first in this direction, even though they are not directly concerned with this problem.  Hasse starts with a number field $k$ and a number of its prime ideals $\mathfrak{p}_i$.  He then finds infinitely many extension fields $K$, in which the $\mathfrak{p}_i$ decompose into prime ideals of prescribed order and multiplicity.  He then significantly tightens this result by taking the prescribed decompositions to be \textit{regular} and requires $K/k$ to be relatively abelian.  He provides the existence proof under certain restrictive conditions.
    \item The general form of problem I can be formulated for relative fields in the following sense (see \cite{ref84}):
    \\ \textit{Problem A.}  Let an algebraic number field $k$ be given, whose [Galois] group [over $\mathbb{Q}$] is given by $\mathfrak{g}$.  Furthermore, let an abstract group $\mathfrak{G}$ and a normal subgroup $\mathfrak{H}$ be given, such that the factor group $\mathfrak{G}/\mathfrak{H}$ is isomorphic to $\mathfrak{g}$.  One must find neccesary and sufficient conditions such that there exists an extension field $K$ of $k$ whose [Galois] group [over $\mathbb{Q}$] is isomorphic to $\mathfrak{G}$.\\
    The following example shows that this problem is not solvable in a few cases.  Let $k$ be cyclic of prime degree $l$.  Let $\mathfrak{G}$ be cyclic of order $l^2$.  The \textit{critical} prime numbers of $k$ are clearly $\equiv 1 (\text{mod } l)$, but not $\equiv 1 (\text{mod }l^2)$.  It follows from the \textit{number theoretic monodromy theorem} that $K$ contains at least one inertia group of order $l^2$ (see \cite{ref82}).  The prime number $p$ corresponding to this inertia group must also be critical in $k$, which is impossible since it is known that it satisfies the congruence $p \equiv 1 \quad (\text{mod }l^2)$.\\
    This example also shows that the solvability of problem A is not just determined by the structure of the group $\mathfrak{G}$ but also via certain arithmetic properties of the field $k$.
    \item Especially important results regarding solving problem A concerning abelian groups were obtained by A. Scholz \cite{ref63}\cite{ref64}.  His investigations are mostly concerned with two-step solvable groups (that is, groups whose commutator subgroups are abelian) and are far-reaching in two directions.  First of all he established a very utilitarian classification of two-tiered groups.  The easiest of his classes, which he called \textit{disposition groups}, allows for a solution to problem A independently of the arithmetic properties of the field $k$ \cite{ref63}.  One can define the disposition group as a group $\mathfrak{G}$, whose normal abelian subgroups $\mathfrak{H}$ are the direct product of all cyclic groups that that are conjugated by them.  In addition each of these cyclic groups is required to have exactly $\mathfrak{H}$ as its normalizer.\\
    Scholz proved the following two theorems on disposition groups:
    \begin{enumerate}[label = \roman*)]
        \item A disposition group $\mathfrak{G}$ is fully determined as soon as one knows the group $\mathfrak{G}/\mathfrak{H}$ and the order of a generating element of $\mathfrak{H}$.
        \item Given an algebraic number field with [Galois] group $\mathfrak{G}/\mathfrak{H}$, one can always find a extension field $K$, whose [Galois] group is isomorphic to $\mathfrak{G}$.
    \end{enumerate}
    Scholz later showed that the first of these  theorems can be expanded to the case where neither $\mathfrak{G}/\mathfrak{H}$ nor $\mathfrak{H}$ are required to be abelian, but are completely arbitrary.\\
    Secondly, Scholz investigated other types of two-tiered groups extensively \cite{ref64}.  He namely found two \textit{maximal types} among all two-tiered groups (that is those types, such that every two-tiered group can be represented as a factor group or a product of groups of this type): \textit{ring groups} and \textit{branch groups}.  Ring groups are always factor groups of certain powers of disposition groups.  Branch groups on the other hand possess a feature that does not allow such a representation: their commutator subgroup is not properly contained in any abelian subgroup.  There the question of whether relative fields with [Galois] branch groups exist remains open.
    \item I have made a contribution towards problem A \cite{ref84}. In this I somewhat broadened the definition of the disposition group by throwing out the requirement that each generating cyclic subgroup of $\mathfrak{H}$ permits no normalizer besides $\mathfrak{G}$.  From this arise the so-called \textit{Scholz groups}, corresponding to the \textit{purely branched} fields (that is those containing no unbranched extension fields of $k$) with relative discriminants, whose prime ideal divisors within $k$ are not critical. The question of the uniqueness of a Scholz group given the factor group $\mathfrak{G}/\mathfrak{H}$, the order of a generating element of $\mathfrak{H}$ and its normalizer, remains open.  The question of the existence of a relative field with a given Scholz group can be reduced to to the existence of principal prime ideals $\mathfrak{p}$ with prescribed values of a Hasse norm remainder symbol$\left(\frac{p, K}{\mathfrak{p}}\right)$ (see Hasse \cite{ref32} III).  However this question exceeds the scope of the currently known analytical ideal theory.
    \item If the group $\mathfrak{G}$ does not belong to the type of Scholz group, then one can only expect the existence of a corresponding extension field when it contains an absolute \textit{Teilklassenk{\"o}rper} or its relatively critical prime ideals are also critical within $k$.  In this case we are subject to the full control of the individual oddities of the field $k$.\\
    The [Galois] groups the of absolute class fields are also not allowed to be completely arbitrary.  On the one hand they are restricted by the aforementioned monodromy theorem, due to which all interia groups generate the full Galois group of the field, while the interia groups of a relatively unbranched field induce an unambiguous mapping on the inertia groups of the base field $k$.  Operating on this principle, I carried out a classification of the possible types of [Galois] groups of absolute class fields \cite{ref82}.  On the other hand, F. Pollaczek \cite{ref60} and Scholz \cite{ref66} discovered and developed a restriction of the absolute class fields based on the properties of the \textit{Grundeinheitensystem}.
\end{enumerate}
\renewcommand{\thesection}{\S\arabic{section}}
\section{On the analytic form of prime numbers belonging to a prescribed permutation class}
\renewcommand{\thesection}{\arabic{section}}
\begin{enumerate}[label = \arabic*.]
    \item It is known that an algebraic number field is not fully determined by its Galois group.  The known invariants that fully determine a field are the so-called Artin symbols $\left(\frac{K}{\mathfrak{p}}\right)$ (\cite{ref32} III, p. 6), that is the permutation classes containing individual prime numbers.  There are infinitely many such invariants, which implies that they cannot be independent from one another.  The existing relations between them can be derived provided we know the analytical form in which we can represent the prime numbers belonging to certain permutation classes, that is corresponding to the same value $\left(\frac{K}{\mathfrak{p}}\right)$.  It can be proven, that such analytical forms depend on the structural properties of the corresponding [Galois] groups $\mathfrak{G}$.  This connection yields a deep insight into the arithmetic structure of a group with known Galois group.
    \item First I will recall the classical case of an abelian field.  For a prime number to belong to a given permutation (in this case each permutation class consists of only one permutation), it must be representable in the form of one of the arithmetic progressions $ax + b$ unambiguously corresponding to the element (permutation) of the Galois group in question.  The number $a$ is one and the same for all permutations and consists of the prime numbers that come up in the driscriminant of the group, while the $b$ are assigned to the different elements (permutations) of the Galois group.
    \item The other, now classic, case is that of complex multiplication of elliptic functions. Given a number field that is relatively abelian over an imaginary quadratic field, then each prime number belonging to a given permutation class can be represented by one or more quadratic forms, whose discriminants depend on the field (more precisely: equal to the discriminant of the imaginary quadratic field) and whose classes are assigned to the permutation classes.
    \item This fact was suspected by Kronecker (\textit{``Jugendtraum''}, or ``youthhood dream'') and proven by R. Fueter \cite{ref24}\cite{ref25}.  The principles that it is based on follow from the general class field theory. Consider only those numbers of a field $k$ which satisfy certain congruences modulo an ideal $\mathfrak{f}$ (which one calls leader\footnote{\textit{F{\"u}hrer} in the original German.-\textit{YS}}.  Taking these numbers as principal ideals, and letting $h$ equal the number of classes in the newly defined sense, it follows that there exists a relatively abelian field $K$ of relative degree $h$.  This field, a so-called class field, has the property a prime ideal of $k$ decomposes fully in $K$ if and only if it lies in the principal class (Ph. Furtw{\"a}ngler, \cite{ref26}; T. Takagi, \cite{ref76}). Conversely we can regard any relatively abelian field over $k$ as a class field with a suitably chosen leader (Fueter, \cite{ref24}; Takagi, \cite{ref76}; Hasse, \cite{ref32}).
    \item We can ask ourselves the general question concerning the analytic form of prime numbers belonging to classes generating by powers of a permutation $S$.  Let a normal algebraic number field $K$ be given, and denote its [Galois] group by $\mathfrak{G}$.  Furthermore let $\mathcal{S}$ be a permutation in $\mathfrak{G}$.  For a prime number $p$ to belong to the classes generated by powers of $\mathcal{S}$, it is necessary and sufficient for the subfield $K_{\mathcal{S}}$ of $K$ belonging to $\mathfrak{Z}_{\mathcal{S}}$ to contain a prime ideal divisor $\mathfrak{p}$ of $p$ of the first degree, where $\mathfrak{Z}_{\mathcal{S}}$ denotes the cyclic subgroup of $\mathfrak{G}$ generating by the powers of $\mathcal{S}$.\\
    Let $\mathfrak{a}_1,\mathfrak{a}_2,...,\mathfrak{a}_k$ be a system of representatives of all distinct ideal classes of $K_{\mathcal{S}}$, and let $(\mu_1^{(i)},\mu_2^{(i)},...,\mu_n^{(i)})$ be respective bases of the ideals $\mathfrak{a}_i (i = 1,2,...,h)$. Then we have that $N(\mu_1^{(i)}x_1 + \mu_2^{(i)z}x_2 + ... + \mu_n^{(i)}x_n)$ is the form (admitting a decomposition) of $n$-th degree in the variables $x_1,x_2,...,x_n$, whose coefficients have the number $N(\mathfrak{a}_i)$ as their greatest common divisor.  The quotient
    \begin{equation}\label{3.1}
        \frac{N(\mu_1^{(i)}x_1 + \mu_2^{(i)}x_2 + ... + \mu_n^{(i)}x_n)}{N(\mathfrak{a}_i)} = f_i(x_1,x_2,...,x_n) \quad (i = 1,2,...,h)
    \end{equation}
    is thus a primitive form of $n$-th degree.  If we allow the $x_i$ to take on any integer values, then the norm $f_i(x_1,x_2,...,x_n)$ takes on the values of the norms of all ideals, whose classes are opposite to the class of $\mathfrak{a}_i$.  Since if $b$ lies in the opposite class to that of $\mathfrak{a}_i$, it follows that $b\mathfrak{a}_i$ is a principal ideal, which may be associated with a number $\mu_1^{(i)}x_1 + \mu_2^{(i)}x_2 + ... + \mu_n^{(i)}x_n$ from the ideal $\mathfrak{a}_i$.  Thus, it holds that
    \[N(b)N(\mathfrak{a}_i) = N(\mu_1^{(i)}x_1 + \mu_2^{(i)}x_2 + ... + \mu_n^{(i)}x_n) = N(\mathfrak{a}_i)f(x_1,x_2,...,x_n).\]
    Given $\mathfrak{p}$ a prime ideal of $K_s$ of first degree, it holds thus that $N(\mathfrak{p}) = p$, one wishes to find the representative $\mathfrak{a}_i$, whose class is opposite to the class of $\mathfrak{p}$. Thus $\mathfrak{p}$ can be represented in the form $f_i(x_1,x_2,...,x_n)$.  Conversely, if $\mathfrak{p}$ can be represented in the form $f_i(x_1,x_2,...,x_n)$, then $N(\mathfrak{a}_i)p$ can be represented in the form $N(\mu_1^{(i)}x_1 + \mu_2^{(i)}x_2 + ... + \mu_n^{(i)}x_n)$. The number $\mu_1^{(i)}x_1 + \mu_2^{(i)}x_2 + ... + \mu_n^{(i)}x_n$ is divisible by $\mathfrak{a}_i$, and the norm of the quotient ideal $\frac{\mu_1^{(i)}x_1 + \mu_2^{(i)}x_2 + ... + \mu_n^{(i)}x_n}{\mathfrak{a}_i}$ is equal to $p$.  It follows from this that there exists an ideal with norm $p$.  This ideal must be a prime ideal of degree 1.
    \item To set up the conditions for a prime number $p$ to be a member of the \textit{division} of $\mathcal{S}$, we must exclude its membership to the powers $\mathcal{S}^k$ of $\mathcal{S}$, whose exponents $k$ are not relatively prime to the order $f$ of $\mathcal{S}$.  For this to occur $p$ in $K_{\mathcal{S}}$ must contain at least one prime ideal of first degree, while this does not happen for any $K_{\mathcal{S}^k}$, $(k,f) \neq 1$.  Given
    \begin{equation}\label{3.2}
        g_i(x_1,x_2,...,x_n),\quad h_i(x_1,x_2,...,x_n),...
    \end{equation}
    the corresponding forms for the fields $K_{\mathcal{S}^k}$, constructed in the same way as $f_i(x_1,x_2,...,x_n)$, we have that $p$ belongs to the division of $\mathcal{S}$ if and only if it admits a representation in terms of one of the forms (\ref{3.1}), but none in terms of the forms (\ref{3.2}).
    \item How can one characterize the membership of a prime number $p$ to the \textit{class} of $\mathcal{S}$?  I can do this only if I know a number $a$ such that $p^f \equiv 1 \text{ (mod $a$)}$, but that this congruence fails for any power of $p$ whose exponent is less than $f$.  Using this to construct the subfield $K(\eta)$ of the $a$-th roots of unity, it follows that $p$ remains irreducible in $K(\eta)$.  If $p\equiv b \text{ (mod $a$)}$, then $\eta^p \equiv \eta^b \text{ (mod $p$)}$  holds. Now we form the element
    \[\xi = \eta\omega  + \eta^b\omega^{\mathcal{S}} + ... + \eta^{b^{f-1}}\omega^{S^{f-1}},\]
    where $\omega$ denotes an element of $K$, and $\Phi(\xi) = 0$ is an equation satisfied by $\xi$.  If one takes as established that $p$ belongs to the \textit{division} of $\mathcal{S}$, then $p$ belongs to the \textit{class} of $\mathcal{S}$ if and only if the congruence $\Phi(\xi) \equiv 0 \text{ (mod $p$)}$ has rational roots.  That is, if $p$ has at least one prime ideal divisor of first degree in $K(\xi)$.  One can hence set up a system of forms such that $p$ can be represented by it if and only if $p$ belongs to the class of $\mathcal{S}$.  One can certainly find the corresponding number $a$ for each $p$; there are however different forms corresponding to different $p$.  I can therefore not set up a finite number of forms that are valid for all prime numbers.
    \item The classically regarded criteria in No. 2 and 3 do not follow from this criterion.  To set up a more general criterion we regard the case where $\mathfrak{G}$ possesses an abelian normal subgroup $\mathfrak{H}$.  Then one can understand $K$ as a relatively abelian field with respect to the field $k$ corresponding to $\mathfrak{H}$.  $K$ is thus a class field of $k$.  By the general reciprocity law of E. Artin \cite{ref2} there is a one-to-one correspondence between the permutations in $\mathfrak{H}$ and the ideal classes (more specifically: the residue classes of a certain ideal class subgroup) of $k$, which has the character of an isomorphism.  The forms (\ref{3.1}) corresponding to $k$ decompose thusly in the system of forms $\mathfrak{B}_i$, of which each represents  one of the aforementioned residue classes.  The number of these systems of forms is equal to the order of $\mathfrak{H}$.  The Artin reciprocity law says that a prime number $p$ can be represented by one of the forms in the system $\mathfrak{B}_i$ if and only if it belongs to the class of $S_i$,  where $S_i$ is one of the permutations in $\mathfrak{H}$ corresponding to the system $\mathfrak{B}_i$.  This system of forms has the advantage that the degree of its forms in generally smaller.  For instance, if $K$ is absolutely abelian, then $k$ is the rational field, so that the norms coincide with the numbers themselves.  The partitioning of numbers of the rational field in the ``narrower'' sense is nothing else but their distribution among the congruence classes by a certain relation, which one calls leader.  If $k$ is quadratic, we arrive at the quadratic forms, in complete concurrence with the general theory.
    \item  We still note that the forms corresponding to the field $k$ allow for the so-called form composition.  For instance, given
    \[N(\mathfrak{a}) = f_1(x_1,x_2,...,x_n), \quad N(\mathfrak{b) = f_2(y_1,y_2,...,y_n)}\]
    then it follows that
    \[N(\mathfrak{a}\mathfrak{b}) = f_1(x_1,x_2,...,x_n) c\dot f_2(y_1,y_2,...,y_n).\]
    On the other hand, if $\mathfrak{a}\mathfrak{b}$ lies in the opposing ideal class to $\mathfrak{a}_3$, then $N(\mathfrak{a}\mathfrak{b}) = f_3(z_1,z_2,...,z_n)$, where $z_1,z_2,...,z_n$ are certain integer-valued bilinear expressions in $x_1,x_2,...,x_n$ and $y_1,y_2,...,y_n$.  One can obtain these by considering the bilinear expression $\sum\limits_{i,j}\mu_i^{(1)}\mu_j^{(2)}x_ix_j$, where $(\mu_1^{(1)},\mu_2^{(1)},...,\mu_n^{(1)})$, $(\mu_1^{(2)},\mu_2^{(2)},...,\mu_n^{(2)})$ denote the bases of the ideals $\mathfrak{a}_1,\mathfrak{b_2}$ respectively.  By expressing the $\mu_i^{(1)}\mu_j^{(2)}$ using a basis $(\mu_1^{(3)},\mu_2^{(3)},...,\mu_n^{(3)})$ of the ideal $\mathfrak{a}_1\mathfrak{b}_2$: $\mu_i^{(1)}\mu_j^{(2)} = \sum\limits_s c^s_{ij}\mu_s^{(3)}$, i.e. $\mathfrak{a}_1\mathfrak{a}_2\mathfrak{a}\mathfrak{b} = \sum\limits_{i,j,s} c_{ij}^sx_ix_j\mu_s^{(3)}$, one then sets the $z_s$ equal to the coefficients of $\mu_s^{(3)}$: $z_s = \sum\limits_{i,j}c_{ij}^sx_ix_j$.  It is easy to understand that this composition of forms corresponds to the multiplication of their corresponding forms.
    \item To obtain a simplest analytic expression of prime numbers belonging to different permutation classes of a field we must find maximal abelian normal subgroups of its [Galois] group $\mathfrak{G}$.  We note that a permutation $\mathcal{S}$ in $\mathfrak{G}$ can be contained in a abelian normal subgroup of $\mathfrak{G}$ if and only if \textit{its class is abelian},  i.e. all permutations in its class commute with each other.  If $\mathfrak{C}_1,\mathfrak{C}_2,...,\mathfrak{C}_k$ are all the abelian classes in $\mathfrak{G}$, we have that each abelian normal subgroup consists of those permutations in the classes which commute with one another.  For instance, given $\mathfrak{C}_1$ and $\mathfrak{C}_2$ commuting, then the \textit{join}, that is the smallest group containing $\mathfrak{C}_1$ and $\mathfrak{C}_2$, is an abelian normal subgroup of $\mathfrak{G}$.  The uniqueness of a maximal abelian normal subgroup cannot be guaranteed since the commutativity of classes is not a transitive property.
    \\
    The following example shows that cases exist where $\mathfrak{G}$ contains multiple different maximal abelian subgroups.  Define $\mathfrak{G}$ using 3 generating elements $s_1,s_2,s_3$, where are subject to the following relations:
    \[s_1^p = s_2^p = s_3^p = 1, s_1s_2 = s_2s_1, s_1s_3 = s_3s_1, s_2s_3 = s_3s_2s_1\]
    ($p$ is a prime number).  Both subgroups $(s_1,s_2)$ and $(s_1,s_3)$ are abelian normal subgroups of $\mathfrak{G}$.  Both are maximal since they are only properly contained in $\mathfrak{G}$, which is not abelian.  On the other hand, they are distinct from one another.
    \item As an example we will consider the general cubic number field $K$.  Its alternating group possesses a quadratic subfield $k$, and by using Fueter-Takagi theory one can regard $K$ as a ring class field of $k$.  Corresponding to the ring classes of $k$ in consideration we have a system of binary quadratic forms, which decompose into the 3 subsystems $\mathfrak{B}_1,  \mathfrak{B}_2, \mathfrak{B}_3$, to which one can associate each of the three permutations in $\mathfrak{H}$.  Given $D$ the discriminant of the system of forms, we have that $\left(\frac{D}{p}\right) = \pm 1$ is the necessary and sufficient condition for $p$ to admit a representation through one of these forms.  The \textit{principal system} is the system from $\mathfrak{B}_1,\mathfrak{B}_2,\mathfrak{B}_3$, say $\mathfrak{B}_1$, which satisfies the property $\mathfrak{B}_1\mathfrak{B}_1 = \mathfrak{B}_1$.  Then all prime numbers relatively prime to the discriminant of $K$ belong to the following three types:
    \begin{enumerate}[label = \roman*)]
        \item $\left(\frac{D}{p}\right) = -1$, $p$ does not belong to $\mathfrak{H}$, hence belonging to one of the transpositions.
        \item $\left(\frac{D}{p}\right) = +1$, $p$ can be represented through one of the forms $\mathfrak{B}_2, \mathfrak{B}_3$.  $p$ belongs to one of the 3-cycles.
        \item $\left(\frac{D}{p}\right) = +1$, $p$ can be represented through one of the forms $\mathfrak{B}_1$.  $p$ belongs to one of the identity permutation.
    \end{enumerate}
    The cubic field was investigated in this respect by Dedekind \cite{ref17}, Voronoi \cite{ref85} and Takagi.  More recently Hasse \cite{ref33} investigated the cubic field from a class field theoretic perspective, by also considering the critical prime ideals.  This work was kindly brought to my attention by B. Delaunay.
    \item I want to mention a very elegant procedure due to A. Speiser \cite{ref73} which serves the purpose of determining the order $\mathfrak{f}$ of a permutation whose class contains a prime number $p$.  I allow myself to modify the proof somewhat.  Given
    \begin{equation}\label{3.3}
        f(x) = x^n + a_1x^{n-1} + a_2x^{n-2} + ... + a_{n-1}x + a_n = 0
    \end{equation}
    the equation to be examined, we consider the difference equation
    \begin{equation}\label{3.4}
        y(m+n) + a_1y(m+n - 1) + ... + a_{n-1}y(m+1) + a_ny(m) = 0.
    \end{equation}
    We know its general solution to be given by
    \begin{equation}\label{3.5}
        y(m) = C_1\alpha_1^m + C_2\alpha_2^m  + ... + C_n\alpha_n^m,
    \end{equation}
    where $\alpha_1,\alpha_2,...,\alpha_n$ are the roots of equation (\ref{3.3}) and $C_1,C_2,...,C_n$ are arbitrary constants.  Setting $y(1) = y(2) = y(n-1) = 0, y(n) = 1$, we have $C_i = \frac{1}{f'(\alpha_i)}$.  Now take $GF[p^{\mathfrak{f}}]$ as a base field.  If $u$ is the smallest (integer-valued) period of the function $y(m)$ modulo $p$, it follows that $\alpha_i^u \equiv 1 \text{ (mod $p$) ($i = 1,2,...,n$)}$, and vice-versa.  This is true since if $y(u) \equiv y(u+1) \equiv.... \equiv y(u + n - 1) \equiv 0 \text{ (mod $p$)}$, $y(u + n) \equiv 1 \text{ (mod $p$)}$, it follows that $C_i \equiv \frac{1}{\alpha_i^uf'(\alpha_i)} \text{ (mod $p$)}$, i.e. $a_i^u \equiv 1 \text{ (mod $p$)}$.  However, since $(\alpha_i \to \alpha_i^p)$ is a generating permutation for the [Galois] group of the congruence $f(x) \equiv 0 \text{ (mod $p$)}$ and its order is hence equal to the smallest number $\mathfrak{f}$, for which $a_i^{p^{\mathfrak{f}}} \equiv \alpha_i \text{ (mod $p$)}$ holds, it follows that $\mathfrak{f}$ is the smallest exponent for which $p^{\mathfrak{f}} \equiv 1 \text{ (mod $u$)}$ holds.
    \item Hasse has recently layed new groundwork for the question concerning the arthimetic structure of number fields, by connecting the theory of number fields with the so-called \textit{``algebras''} (i.e. \textit{hypercomplex systems}) \cite{ref34}\cite{ref9}.  Each field has a corresponding \textit{``cyclic''} hypercomplex system, introduced by E. Noether under the name \textit{crossed product}.  However, we conversely have that multiple fields with different Galois groups, in particular cyclic fields, correspond to the same cyclic algebra, meaning this is the deepest connection in which we can bring arithmetic properties of fields those of cylic fields.
\end{enumerate}
\renewcommand{\thesection}{\S\arabic{section}}
\section{The resolvent problem}
\renewcommand{\thesection}{\arabic{section}}
\begin{enumerate}[label = \arabic*.]
    \item There is a question in the general algebraic field theory which includes the resolvent problem as a special case:\\
    Let $K$ be a given field of a rational functions in the variables $x_1,x_2,...,x_n$.  One must determine the \textit{true transcendence degree} of $K$ in relation to a certain one of its subfields $k$.  That is, one finds the smallest number $m$ such that $K$ can be expressed as a direct product of a field of \textit{algebraic} functions in $m$ variables, which split in $k$, and a certain subfield of $k$ (see below \S5).\\
    To clarify the connection between this question and the resolvent problem, we will consider the the field $K$ of all rational functions in the variables $x_1,x_2,...,x_n$, while $k$ consists of the elementary symmetric functions $a_1,a_2,...,a_n$ of the variables $x_1,x_2,...,x_n$ and their rational functions.  One must find an equation of $n$-th degree \textit{(resolvent)} whose coefficients are rational functions of $a_1,a_2,...,a_n$ and as few as possible ($m$) are functionally independent, and whose roots generate the whole field $K$.  In other words, we are dealing with the Tschirnhaus transformation of the equation
    \begin{equation}\label{4.1}
        x^n + a_{1}x^{n-1} + ... + a_{n-1}x + a_n = 0
    \end{equation}
    with the unrestricted variable coefficients $a_1,a_2,...,a_n$ into an equation, whose coefficients possess a number of degrees of freedom that are as small as possible.\\
    It is of use to expand this task by introducing irrational functions of $a_1,a_2,...,a_n$ into the \textit{coefficient field} $k$, however are also not of a higher transcendence degree than $m$.  An example of this is $\sqrt{D}$, where $D$ denotes the discriminant of equation (\ref{4.1}).\\
    This problem can be summarized as a natural extension of the original task of Galois theory, the solution by radicals.  This is because the representation of roots through radical expressions has the advantage that it permits the determination of roots through a sequence of operations that each only is concerned with one variable.  To this end one can create a table that assigns to each radicand its radical, such that these tables enable the computation of roots of all solvable equations of a given degree.
    \item The aforementioned property is in no way characteristic of solvable equations.  One can rather set up the detection of roots through operations of a totally different kind, for which this property applies to all of them.  First of all one can mention the general equation of fifth degree.   It was known for a long time (Bring, \cite{ref12}), that one can convert it into the so called Bring-Jerrard form
    \begin{equation}\label{4.2}
        y^5 + py + q = 0
    \end{equation}
    by applying a Tschirnhaus transformation, whose coefficients are roots of the equations of fourth degree, that is permit representations via radical expressions (see J. J. Sylvester, \cite{ref75}).  On the other the hand the connection between equations of fifth degree and the division problem of periods of elliptic functions (thus the arguments of modular functions) into 5 is well known.  This allows the general equation of fifth degree to be solved in a transcendental way (see \cite{ref29},\cite{ref86}).
    \item This problem was solved by F. Klein in a way that allows a perspective into the general resolvent problem \cite{ref41}.  He namely converted the general equation of fifth degree into a slightly different normal form
    \begin{equation}\label{4.3}
        y^5 + 15y^4 - 10\gamma\cdot y^2 + 3\gamma^2 = 0
    \end{equation}
    via the use of only quadratic irrationalities, one of which is $\sqrt{D}$ and the other is $\sqrt{5}$.\\
    The second, much more important contribution from Klein towards the resolvent problem consists of developing the deeper reason for why the resolvent problem is solvable for equations of fifth degree.  Specifically he linked this problem to the so-called \textit{form problem}, which consists of the following.  If we regard the largest finite group $\mathfrak{G}$ of binary linear transformations, the icosahedral group, one can associate to it an equation of degree 60
    \begin{equation}\label{4.4}
        (\mathfrak{D}^{30} + 522\mathfrak{D}^{25} - 10005\mathfrak{D}^{20} - 10005\mathfrak{D}^{10} - 522\mathfrak{D}^{5} + 1)^2 = z\cdot \mathfrak{D}^5(\mathfrak{D}^{10} + 11\mathfrak{D}^5 - 1)^5 
    \end{equation}
    whose coefficients depend on a form $z$ of the variables $x_1,x_2$, which are invariant under the transformations of $\mathfrak{G}$, while the roots are mapped to each other via these transformations.  The Galois group of this equation is isomorphic as the icosahedral group to the alternating group of fifth degree.  From this one can conclude that each equation of fifth degree can be converted to the form (\ref{4.4}) (or also (\ref{4.3})) by using a rational transformations, whose coefficients may contain $\sqrt{D}$, where $D$ denotes the discriminant of this equation.\\
    The basic idea of reducing fifth degree equations to a parametric resolvent comes from the fact that the composition series of the symmetric group of fifth degree consists of two factors, where one is a group of degree 2, while the other is isomorphic to the icosahedral group, which can be represented as a group of linear fractional transformations.
    \item This idea was applied by Klein to other equations, in particular the simple group of order 168, which can be represented using ternary linear homogeneous transformations.  This group corresponds to a special class of equations of 7-th degree that have this group as their Galois group.  Since the group of ternary linear homogeneous transformations is isomorphic to the group of linear fractional transformations of two variables, we can similarly conclude that the equations we're concerned with possess a two parameter resolvent.\\
    Slightly more complicated was the matter of the alternating equations of 6-th degree.  The alternating group of 6-th degree does not possess a representation by ternary linear homogeneous transformations.  It was for this reason that Klein suspected that the general equation of 6-th degree does not possess a 2 parameter resolvent, and sought out 3 parameter resolvents.  A. Wiman \cite{ref89} noticed, that in spite of that this group admits a representation in terms of linear fractional transformations in two variables.  This is because one can construct this group as a factor group of a certain group of order 1080, which can be represented as a linear homogeneous group of three variables.  Its corresponding normal subgroup $\mathfrak{H}$ of order 3 lies in the center of the group and thus appears as a group whose transformations only affect the multiplication of the variables with the 3rd roots of unity.  Now looking at how this affects the relations between the variables, we have that they form a group of linear fractional transformations.  Acting on the original variables with the transformations in $\mathfrak{H}$, we see that the relations remain unchanged.  The group of linear fractional transformations we have just constructed it thus isomorphic to the factor group, that is the altnerating group of 6-th degree.
    \item The question regarding the representation of finite groups by linear fractional transformations was investigated in general by I. Schur \cite{ref69}\cite{ref70}.  It was proven that this task could always be completed using a finite number of operations.  That is, in order to find all such representations of a given finite group $\mathfrak{G}$ one must find a covering group  $\mathfrak{K}$ corresponding to $\mathfrak{G}$ possessing the following properties:
    \begin{enumerate}[label = \Roman*.]
        \item $\mathfrak{G}$ is isomorphic to a factor group $\mathfrak{K}/\mathfrak{M}$.
        \item $\mathfrak{M}$ is contained in the center of $\mathfrak{K}$.
        \item The commutator subgroup of $\mathfrak{G}$ is the factor group $\mathfrak{D}/\mathfrak{M}$, $\mathfrak{D}$ denoting the commutator subgroup of $\mathfrak{K}$.
    \end{enumerate}
    Each group $\mathfrak{G}$ corresponds to only a finite number of distinct covering groups $\mathfrak{K}$.  If one finds all representations of $\mathfrak{K}$ by linear homogeneous transformations, each such representation induces a representation of $\mathfrak{G}$ by linear fractional transformations.  This is the general method for the formation of all representations of $\mathfrak{G}$ of this kind.
    \item A. Wiman \cite{ref90} studied the question of representing the symmetric and alternating groups of higher degrees.  It was shown that for $n \geq 8$ the symmetric groups $\mathcal{S}_n$ can be represented as linear homogeneous groups of no less than $n-1$ variables, with the same statement holding true for the alternating groups.
    \item D. Hilbert \cite{ref36} posed a problem (``13. problem'') that has many points of contact with the Klein problem discussed above.  Let the roots of equation (\ref{4.1}) be regarded as functions of $n$ variables.  One asks if these functions can be written as superpositions of functions of a smaller number $k$ of variables and rational operations.  Later \cite{ref37} he found the following values for $k$ for $n \leq 9$:  
    \[
    \begin{tabularx}{\textwidth}{ X|X|X|X|X|X| }
    $n$ & 5 & 6 & 7 & 8 & 9 \\
    \hline 
    $k\leq$ & 1 & 2 & 3 & 4 & 4  \\
    \end{tabularx}.
    \]
    Wiman \cite{ref91} generalized this result by showing that for all $n\geq 9$ we have the inequality $k \leq n-5$, i.e. any general equation of degree $n \geq 9$ can be reduced by at least 5 variables via application of a Tschirnhaus transformation. (See also R. Garver, \cite{ref28}).
    \item I \cite{ref83} have set myself the goal to examine the connection between the reducibility of an equation and the ability to represent their Galois group as a transformation group of as few variables as possible somewhat closer.  To this I will introduce the term \textit{Einkleidungsgruppe}\footnote{This translates as \textit{outfitting group} or \textit{investiture group.}-\textit{YS}} (short E.G.):
    \\
    Given $\mathfrak{G}$ a finite group, we call a continuous group $\Gamma$ an E.G. of $\mathfrak{G}$ if and only if it satisfies the following conditions:
    \begin{enumerate}[label = \roman*)]
        \item $\Gamma$ contains a subgroup isomorphic to $\mathfrak{G}$.
        \item There is no proper subgroup of $\Gamma$ possessing property 1).
        \item There is no proper factor group of $\Gamma$ possessing property 1).
    \end{enumerate}
    I then proved the following \\
    \textit{Theorem.}  An algebraic equation with unrestricted variable coefficients possesses a $k$ parameter resolvent if and only if its Galois group has an E.G. that can be represented in $k$-dimensional space. \\
    To prove the first part of this theorem we will assume that the group $\mathfrak{G}$ is simple, which is sufficient for our purposes.  In this case it follows that each E.G. is also a simple group.  Since if a continuous group $\Gamma$ containing $\mathfrak{G}$ as a subgroup contains a proper normal subgroup $\Gamma_1$, it follows that $\Gamma_1$ either contains all of $\mathfrak{G}$ or not elements of $\mathfrak{G}$ besides the identity.  In the first case this contradicts condition 2).  In the second case it follows that $\Gamma/\Gamma_1$ has a subgroup isomorphic to $\mathfrak{G}$, contradicting condition 3).
    \\
    One can moreover show the following: If a continuous group $\Gamma$ having $\mathfrak{G}$ as a subgroup can be represented in $k$-dimensional space (abbreviated: $\Gamma$ is a $k$-group), then the subgroup of $\Gamma$ which is an E.G. of $\mathfrak{G}$ is also a $k$-group.  It is evident that a subgroup of a $k$-group is a $k$-group.  To show that a simple factor group of a $k$-group is once again a $k$-group, we recall a theorem initially proven in generality by E.E. Levi \cite{ref44}:
    If a continuous group $\Gamma$ has a simple factor group $\Gamma_1$, it also has a subgroup isomorphic to $\Gamma_1$.
    \item If $\mathfrak{G}$, presented as transformations of the roots as variables, permits an E.G. that is not just \textit{isomorphic}, but \textit{similar} to a $k$-group, one immediately obtains a solution to the problem.  One can namely set up $k$ functions $z_1,z_2,...,z_k$ of the variables $x_1,x_2,...,x_n$ which describe a group isomorphic to $\Gamma$, by applying the transformations in $\Gamma$ to the $x_i$.  In doing so, due to the investigations of Lie-Engel (\cite{ref45}, p. 522), one can find the $z_i$ in a rational way, since they form systems of imprimitivity of $\Gamma$.  If the group $\Gamma$ is transitive, one can take the $z_i$ to be rational functions of $x_1,x_2,...,x_n$.  Otherwise they depend on an irrationality $\theta$ which itself depends on the group $\Gamma$.  If for instance we take $z_i$ as a rational function of $x_1,x_2,...,x_n$, then they generate all elements of the field coinciding with $k[x_1,x_2,...,x_n]$ conjugate to $z_i$.  On the other hand these elements depend solely on $z_1,z_2,...,z_n$, since they are generated by application of the transformations in the group $\mathfrak{G}$, while the transformations in $\mathfrak{G}$ appear in the group $\Gamma$ as transformations taking the $z_i$ into functions of the $z_i$.  It follows from this that $z_i$ satisfies an equation whose coefficients depend on only $k$ parameters, while the $x_i$ can be expressed rationally by the $z_i$.
    \item If $\Gamma$ is an E.G. of the group $\mathfrak{G}$ and isomorphic, but not similar to a $k$-group, one can apply Cartan's principle, which consists of the following \cite{ref14}:
    \\
    If two continuous groups $\Gamma$ and $\Gamma_1$ are isomorphic, one can extend the group $\Gamma$ in such a way that one obtains the group $\Gamma_1$ by applying the transformations of $\Gamma$ to the variables corresponding to the functions of the extended group.
    \\
    By ``extension of a group'' one understands the following.  If $x_1,x_2,...,x_n$ are the variables corresponding to the group $\Gamma$, one introduces new sequences of variables
    \begin{equation}\label{4.5}
        \begin{matrix}
        x_1, & x_2, & ..., & x_n, \\
        x_1^{(1)}, & x_2^{(1)}, & ..., & x_n^{(1)}, \\
        ... & ... & ... & ... \\
        x_1^{(m-1)}, & x_2^{(m-1)}, & ..., & x_n^{(m-1)},
        \end{matrix}
    \end{equation}
    so that one obtains the the extended group $\Gamma$ by simultaneously applying transformations in the original group $\Gamma$ to each sequence:
    \[x_1^{(\lambda)}, x_2^{(\lambda)}, ..., x_n^{(\lambda)} \quad (\lambda = 0,1,2,...,m-1).\]
    \\
    To construct a $k$ parameter resolvent in this case, instead of the variables $x_1^{(\lambda)}, x_2^{(\lambda)}, ..., x_n^{(\lambda)} (\lambda = 0,1,2,...,m-1)$ we introduce the new variables $a_1^{(\lambda)}, a_2^{(\lambda)}, ..., a_n^{(\lambda)}$ using the formulae
    \begin{align}\label{4.6}
        &x_1^{(\lambda)} = a_0^{(\lambda)} + a_{1}^{(\lambda)}x_1 + ... + a_{n-1}^{(\lambda)}x_{1}^{n-1}, \\
        &x_2^{(\lambda)} = a_0^{(\lambda)} + a_{1}^{(\lambda)}x_2 + ... + a_{n-1}^{(\lambda)}x_{2}^{n-1}, \nonumber \\
        &... \quad \quad \quad \quad \quad \quad \quad \quad \quad\quad \quad \quad \quad \quad \quad \quad \quad (\lambda = 1,2,...,m-1)\nonumber\\
        &x_n^{(\lambda)} = a_0^{(\lambda)} + a_{1}^{(\lambda)}x_n + ... + a_{n-1}^{(\lambda)}x_{n}^{n-1}. \nonumber
    \end{align}
    If we apply the transformations $U$ in $\Gamma$ to the variables in (\ref{4.5}), the $a_i^{(\lambda)}$ undergo certain transformations that become the identity when we replace the transformations $U$ with transformations from $\mathfrak{G}$.  Plugging these expressions into one of the functions
    \[z_i(x_1,x_2,...,x_n;x_1^{(1)},x_2^{(1)},...,x_n^{(1)};...;x_1^{(m-1)},x_2^{(m-1)},...x_n^{(m-1)}),\]
    say $z_1$, which is subject to a transformation from a $k$-group whenever one applies a transformation from $\Gamma$ to the $x_i^{(\lambda)}$, one finds that the functions $z_1,z_1^{s_1},...,z_1^{s_{N-1}}$ only depend on $k$ parameters, where $\mathfrak{G} = 1 + s_1 + ... + s_{N-1}$.  Since the $a_i^{(\lambda)}$ remain invariant under these transformations we can replace them with arbitrary rational numbers, provided one requires that the differences $z_1^{s_i} - z_1$ $(i = 1,2,...,N-1)$ are all non-zero.  From this we obtain a Galois resolvent with $k$ parameters.  However, in order to obtain $k$-parameter equation from this that comes from performing a Tschirnhaus transformation on equation $(\ref{4.1})$, we must construct a function from the $z_1^{s_i}$ that belongs to the group that $x_1$ belongs to.
    \item It is not out of the question that the transition of the $x_i$ to the $Z_i$ is not rational, but instead contains an irrationality $\theta$, which itself depends on the invariants of the extended group $\Gamma$.  We will refer to the irreducible equation $R(\theta) = 0$ satisfied by $\theta$ as the \textit{near resolvent}.  This raises the important question of whether or not a near resolvent contains more essential parameters than the equation (\ref{4.1}) itself.  I cannot answer this question in general today.  It is only known that for symmetric $\mathfrak{G}$ the near resolvent has 0 parameters, thus is numerical, and for alternating $\mathfrak{G}$ has at most 1 parameter.
    \item Now we will prove the second part of the stated theorem.  Let $Z_1,Z_2,...,Z_n$ be the roots of a $k$-paramter resolvent of equation (\ref{4.1}), which is written in the fashion of a normal equation.  The $Z_i$ are functions of $x_1,x_2,..,x_n$ belonging to the unity group.  Acting on the $x_i$ using the permutations in the group $\mathfrak{G}$, we obtain induced permutations on the $Z_i$ which form a group isomorphic to $\overline{\mathfrak{G}}$, which differs from $\mathfrak{G}$ only by a different designation of the variables.\\
    At first we will take $Z_1,Z_2,...,Z_n$ as independent variables and will attach to $\overline{\mathfrak{G}}$ the continuous group $\Gamma$ in the following way:  Consider a system $A,B,...,$ a system of permutations in $\overline{\mathfrak{G}}$ which generate the whole group $\overline{\mathfrak{G}}$ using composition.  We consider each of these permutations, say $A$, as a linear homogeneous transformation and depict them in a normal form:
    \[u_i \to \epsilon^{k_i}u_i \quad \quad (i = 1,2,...,n),\]
    where $\epsilon = e^{\tfrac{2\pi i}{m}}$ is a root of unity and $u_i$ denotes a linear function in $Z_1, Z_2,...,Z_n$.  We then take
    \[u_i \to e^{k_it}u_i \quad \quad (i = 1,2,...,n)\]
    as the generating transformation for the group $\Gamma_1$, which then becomes $A$ after setting $t = \frac{2\pi i}{m}$.  Returning to the original variables $Z_1,Z_2,...,Z_n$ and do this with all generating transformations $A,B,...$, we obtain a continuous group $\Gamma$ which is an $E.G.$ for $\overline{\mathfrak{G}}$.  $\Gamma$, being a linear group of the $n$ variables, only has a finite number of parameters.\\
    Next we regard the $Z_i$ as functions of the $x_i$ and view the $x_i$ as coordinates of a space $\mathfrak{R}$, whereby two points $(x_1,x_2,...,x_n)$ and $(x_1',x_2',...,x_n')$ are considered not distinct from one another if and only if
    \[Z_i(x_1,x_2,...,x_n) = Z_i(x_1',x_2',...,x_n') \quad \quad (i = 1,2,..,n).\]
    Since amongst the $Z_1,Z_2,...,Z_n$ there is only a set of $k$ which is functionally independent, it follows that the space $\mathfrak{R}$ only has $k$ dimensions.  The group $\Gamma$ induces a continuous group in $\mathfrak{R}$, that is locally\footnote{See previous footnote concerning isomorphisms \textit{``im kleinen''}. -\textit{YS}} isomorphic to $\Gamma$\cite{ref67}, and contains the group isomorphic to $\mathfrak{G}$ as a subgroup, since for instance $Z_1$ belongs to the identical group of $\mathfrak{G}$.  Thusly we have that $\mathfrak{G}$ admits a $k$-group as an E.G., QED.
    \item In the proof of this theorem it is important that every coordinate system of the parameters of the group $\Gamma$ determines the point of the space $\mathfrak{R}$ in a well-defined manner.  For that reason the Schreier definition \cite{ref67} applies to our group, meaning all of Schreier theory is applicable.  However, if this condition is not satisfied, we could for instance encounter a scenario in which a non-cyclic monodromy group of an algebraic function of a single variable admits a single-parameter continous group as an E.G.  As an example, the continuous group defined by the equations
    \[x + y + z = C_1, \quad x^2 + y^2 + z^2 = C_2, \quad x^3 + y^3 + z^3 = C_3\]
    contains the symmetric permutation group of third degree, which is not even abelian.
    \item Let a finite group $\mathfrak{G}$ without a center be given.  How does one go about determining the E.G. belonging to it that is representable in as few dimensions as possible?  To answer this question one should note that each of the sought after continuous groups is always isomorphic ``im kleinen'' to a linear homogeneous group $\Gamma$.  At the same time we have that $\Gamma$ must either contain $\mathfrak{G}$ itself as a subgroup, or it must have a finite subgroup which contains $\mathfrak{G}$ as the factor group with respect to its center.  The follows from the theory of groups ``isomorphic im kleinen'' due to Schreier \cite{ref67}, according to which all groups ``isomorphic im kleinen'' are factor groups of a covering group with respect to discrete subgroups lying in the center of the covering group.  It then follows from the work of I. Schur \cite{ref69}\cite{ref70} that one must consider all ``covering groups''.  Their quantity is known to be finite.
    \\
    This question admits a fairly straightforward answer in the case where $\mathfrak{G}$ is a simple group.  Since in this case we have that its E.G. is also simple.  However, from the work of W. Killing \cite{ref40} and E. Cartan \cite{ref13}, it follows that besides a finite number of easily specified exceptions there exist only three types of simple groups: 1) full unimodular linear groups ; 2) orthogonal groups; 3) symplectic groups.  It is otherwise known (Cartan, \cite{ref13}) that $n$-dimensional groups of type 1) are $(n-1)$-groups and those of types 2) and 3) are $(n-2)$-groups.  It follows that one only needs to examine linear homogeneous groups of at most $n-1$ dimensions for possible E.G.'s for $\mathfrak{G}$, where $n$ is the degree of (\ref{4.1}).  On the other hand, Wiman (\cite{ref91}; see also R. Garver, \cite{ref28}) proved that the altenerating group of $n$-th degree ($n \geq 8$) cannot be represented by linear homogeneous transformations of degree less than $(n-1)$.  This does not yet allow us to see the task of finding a resolvent in less than $(n-3)$ parameters as impossible.  This is because one must consider the covering groups themselves and not just representability of the alternating group.  We have convinced ourselves using the $n = 6$ case that this sometimes allows a reduction in the number of parameters in the resolvent.  However, I must once again particularly emphasize that the Sylvester-Hilbert-Wiman problem, which for $n \geq 9$ allows a reduction of at least 5 parameters, cannot be regarded as a special case of the Klein problem.  In other words we cannot claim a priori that a reduction in parameters assisted by a chain of resolvents can only be accomplished using a single resolvent.  The theorem valid in classical Galois theory regarding natural irrationalities cannot be directly extended to the resolvent problem.  The accomplishment of the task via chains of resolvents requires an extensive study of the resolvent problem in the case where the coefficients of equation (\ref{4.1}) are not free, but connected by some relations. (See \S 5, No. 8).  In that case it can arise that the finite group to be embedded in an E.G. does not coincide with the monodromy group of equation (\ref{4.1}).  In No. 13 we saw an example of this phenomenon.
    \item The \textit{Einkleidungsproblem} for finite groups by continuous groups is of significant interest in itself.  I cannot say at this point whether a finite group admits a finite or infinite number of non-isomorphic continuous groups as E.G.s.  The representation theory of continuous groups is of use for the solution to this problem (Cartan, \cite{ref13}; Schur, \cite{ref71}; R. Brauer, \cite{ref8}; H. Weyl, \cite{ref88}).  It is however very inconvenient that each continuous group admits infinitely many irreducible linear homogeneous representations.
    \\
    The converse of this problem was already posed a while ago.  This is the problem of classifying all finite groups, that are contianed in a given continuous group.  C. Jordan \cite{ref39} proved the following fundamental theorem for this problem:
    \\
    Each continuous linear homogeneous group $\Gamma$ contains only those non-isomorphic finite subgroups $\mathfrak{G}$ whose factor groups by abelian normal subgroups belong only finitely many times to the group $\Gamma$.
    \end{enumerate}
    \renewcommand{\thesection}{\S\arabic{section}}
    \section{Further questions in the general theory of fields}
    \renewcommand{\thesection}{\arabic{section}}
    \begin{enumerate}[label = \arabic*.]
        \item The questions in the theory of algebraic number fields that regard the algebraic numbers with respect to their rationality are usually solvable using Galois theoretic methods.  However, if the fields in question contain certain transcendental (variable) elements, then the corresponding structural questions for these fields do not admit a straightforward group theoretic description.  Even though I nonetheless include these questions in Galois theoretic body of ideas, I do it for the following reasons: It is first of all not natural to define Galois theory as the body of ideas whose problems can be solved using the Galois group, since the Galois group is a solution tool.  That is to say that these same problems may be solved using significantly different means, meaning a solution tool is in no way suitable for delimiting an area of study.  Secondly we cannot say a priori whether a problem in consideration is not solvable using an appropriately defined notion of the Galois group.  It is far more appropriate to define Galois theory as the body of problems, whose problems concern themselves with the rational dependence of fields and individual field elements.
        \item \textit{Identity of two algebraic fields.}  A field is by no means determined by its Galois group.  On the contrary one can construct distinct fields with isomorphic [Galois] groups.  The question concerning the identity of fields thus actually lies outside of the realm of Galois theory.  If $K$ is a number field this question is essentially number-theoretic.  Its answer would be best achieved if we construct a complete system of invariants for number fields.  For the latter problem there are two known methods.  The first follows from Dedekind-Frobenius theory (see \S3).  This method is inconvenient due to the infinite number of invariants, that are linked to one another using very few known relations.  The second method is based on the behavior of the field discriminant. However, there is no known invariant system of this kind which determines a field uniquely.  Besides this we have that the field discriminant cannot take on every integer value, and to this day we do not know what number values it can take on. 
        \item There is however a purely algebraic method for the solution of the identity problem.  Given $K$ and $K_1$ fields to consider, with isomorphic [Galois] groups $\mathfrak{G}$, we have in general that the compositum $KK_1$ has the direct product $\mathfrak{G}\times\mathfrak{G}$ as its Galois group.  If $K$ and $K_1$ have an irrational intersection, we have that the [Galois] group of $KK_1$ becomes a proper subgroup of $\mathfrak{G}\times\mathfrak{G}$.  In particular, if $K$ and $K_1$ are identical, it follows that the [Galois] group of $KK_1$ is isomorphic to $\mathfrak{G}$.  From this one can derive derive a usesful criterion for whether $K$ and $K_1$ are identical.  Specifically, if 
        \begin{equation}\label{5.1}
            x^n + a_1x^{n-1} + ... + a_{n-1}x + a_n = 0, y^n + b_1y^{n-1} + ... + b_{n-1}y + b_n = 0
        \end{equation}
        are the equations whose individual roots generate the fields $K$ and $K_1$ respectively, it follows that $K = K_1$ if and only if one of the quantities
        \[x_1^ky_1 + x_2^ky_2 + ... + x_n^ky_n \quad \quad (k = 1,2,...,n-1)\]
        is rational, where $x_1,x_2,...,x_n$ and $y_1,y_2,...,y_n$ denote the roots of the equations in (\ref{5.1}).  Each of the quantities $x_1^ky_1 + x_2^ky_2 + ... + x_n^ky_n$ satisfies an equation of degree $n!$, whose coefficients are rational expressions of the $a_i,b_i$.  The fields $K$ and $K_1$ are thus identical if and only if this equation possesses at least one rational root.  Moreover, if $\mathfrak{G}$ has known normal subgroups, one can reduce the degree of this (``mixed'') equation using mixed resolvents \cite{ref78}.
        \item We will examine the $n = 3$ and $n = 4$ cases more closely.  Given
        \[x^3 + px + q = 0, y^3 + \overline{p}y + \overline{q} = 0\]
        as the equations under consideration, we first require that the product of their discriminants is a perfect square.  The quantity $z = x_1y_1 + x_2y_2 + x_3y_3$ satisfies on the mixed equations
        \[z^3 - 3p\overline{p}z - \frac{27}{2}q\overline{q} \pm \sqrt{D\overline{D}} = 0.\]
        Provided that $z$ is one of its rational roots, and if $z^3 - p\overline{p} \neq 0$, one can determine the coefficients of the rational map $y = \alpha_0 + \alpha_1x + \alpha_2x^2$ from the equations
        \[3\alpha_0 - 2p\alpha_2 = 0, -2p\alpha_1 - 3q\alpha_2 = z_1, -2p\alpha_0 -3q\alpha_1 + 2p^2\alpha_2 = u,\]
        where $u  = x_1^2y_1 + x_2^2y_2 + x_3^2y_3 = \frac{3(q\overline{p}z - p^2\overline{q})}{z^2 - p\overline{p}}$. \\
        Now we will consider the case $n = 4$.  Given equations
        \[x^4 + p_2x^2 + p_3x + p_4 = 0, \overline{x}^4 + \overline{p}_2\overline{x}^2 + \overline{p}_3\overline{x} + \overline{p}_4 = 0\]
        to consider, we will initially solve the problem for the cubic equations
        \[z^3 - p_2z^2 - 4p_4z - p_3^2 + 4p_2p_4 = 0, \overline{z}^3 - \overline{p}_2\overline{z}^2 - 4\overline{p}_4\overline{z} - \overline{p}_3^2 + 4\overline{p}_2\overline{p}_4 = 0,\]
        which we know are satisfied by the quantities $z = x_1x_2 + x_3x_4$ and $\overline{z} = \overline{x}_1\overline{x}_2 + \overline{x}_3\overline{x}_4$ respectively.  Next, if we introduce the expressions
        \[\zeta = z_1\overline{z}_1 + z_2\overline{z}_2 + z_3\overline{z}_3,\quad u = z_1^2\overline{z}_1 + z_2^2\overline{z}_2 + z_3^2\overline{z}_3, \quad \overline{u} = z_1\overline{z}_1^2 + z_2\overline{z}_2^2 + z_3\overline{z}_3^2\]
        it follows that the equation
        \[F(T) = T^4 - (2p_2\overline{p}_2 + 2\zeta)T^2 - 8p_3\overline{p}_3T -\tfrac{1}{3}\zeta^2 - \tfrac{8}{3}p_2\overline{u} - \tfrac{8}{3}\overline{p}_2u + \tfrac{14}{8}p_2\overline{p}_2\zeta + p_2^2\overline{p}_2^2 + 16p_2^2\overline{p}_4 + 16\overline{p}_2^2p_4 + \tfrac{64}{3}p_4\overline{p}_4 = 0\]
        has at least one rational root.  If we also have $F'(T) \neq 0$, it follows that one can determine the coefficients $\alpha_0,\alpha_1,\alpha_2,\alpha_3$ of the map $x = \alpha_0 + \alpha_1x + \alpha_2x^2 + \alpha_3x^3$ from the equations
        \[4\alpha_0 - 2p_2\alpha_2 - 3p_3\alpha_3 = 0, -2p_2\alpha_1 - 3p_3\alpha_2 + (2p_2^2 - 4p_4)\alpha_3 = T,\]
        \[-2p_2\alpha_0 - 3p_3\alpha_1 + (2p_2^2 - 4p_4)\alpha_2 + 5p_2p_3\alpha_3= \theta,\]
        \[3p_3\alpha_0 + (2p_2^2 - 4p_4)\alpha_1 + 5p_2p_3\alpha_2 + (-2p_2^3 + 3p_3^2 + 6p_2p_4)\alpha_3 = Z,\]
        where:
        \[\theta = x_1^2\overline{x}_1 + x_2^2\overline{x}_2 + x_3^2\overline{x}_3 + x_1^2\overline{x}_1,\]
        \[Z = x_1^3\overline{x}_1 + x_2^3\overline{x}_2 + x_3^3\overline{x}_3+x_4^3\overline{x}_4.\]
        \item Given $K$ an algebraic function field, we have that Galois theory admits no applications to the solution of the identity problem.  If $K$ possesses just one independent variable, this problem is solvable with the aid of function-theoretic means.  If its genus is namely $p > 1$, it follows that $K$ is determined by $3p - 3$ independent continuous parameter values.  For the equations generating such fields there are known normal forms.  If we find a normal form for each of the fields to be compared, the problem is solved by their compilation (\cite{ref3}, p. 90-92).
        \\
        If $K$ depends on multiple independent variables then the problem remains unsolved in general.  The old German and the Italian geometers obtained several extraordinarily important results in this area.  Yet I fear that even today, regarding this, we must repeat the following phrase, contained in the Z{\"u}rich lecture of F. Enriques \cite{ref19}:\\
        \guillemotleft Malheuresement la plupart de ces probl{\`e}mes demeurent aujourd'hui sans r{\'e}ponse, et les contributions qu'on a port{\'e}es dans ce champ de recherches resemblement en v{\'e}rit{\'e} {\`a} de rares flambeaux au milieu d'une obscurit{\'e} {\'e}paisse.\guillemotright
        \newline
        [\textit{English:} Unfortunately, the majority of these problems remain unanswered, and the contributions we have reached in this area of research truly resemble rare torches in the middle of a thick darkness.]
        \item \textit{Rational minimal basis.} If we consider the field $K_n$ of all rational functions in $n$ independent variables $x_1,x_2,...,x_n$, the question seeking all possible types of its subfields arises.  One can easily set up ``trivial'' examples of such subfields: that is one takes some number $m \leq n$ of functionally independent elements of $K_n$ (i.e. rational functions of $x_1,x_2,...,x_n)$ as generating elements of a subfield.  The field constructed in this manner is obviously either isomorphic to $K_n$ or with the field $K_m$ of rational functions in $m$ ($m < n$) variables.  One can however expect the existence of other kinds of fields:  one can take a subfield generated by some number elements in $K_n$ that may be linked with one another by algebraic relations.  The question arises whether the ``trivial'' instances of subfields comprise all subfields exhaustively.  In other words, one asks after the existence of a system of independent elements of the subfield, by which all elements of this subfield may be represented.  Such a system is called a \textit{rational minimal basis}.\\
        This question was answered affirmingly for $n = 1$ by L{\"u}roth (\cite{ref48}; see also E. Netto, \cite{ref53}).  G. Castelnuovo \cite{ref15} extended this result to case of $n = 2$.  His proof relies on the methods of algebraic geometry.  In the case $n = 3$ we have a counterexample due to G. Fano \cite{ref21} and F. Enriques \cite{ref20}, which they obtained by constructing a subfield of a field of rational functions in three variables admitting no rational minimal basis.
        \\
        The problem of the rational minimal basis has an application in the classical Galois theory, namely to the question after the existence of fields with prescribed Galois groups (\S2).  To answer the latter question it is necessary to solve the problem of the rational minimal basis solely in the case where the subfield under consideration contains the field of elementary symmetric functions of a system of generating elements for $K_n$.  Given this restriction we have that the problem of the rational minimal basis is neither settled nor refuted (see \cite{ref74}).
        \\
        This problem can also be understood in terms of the problem of identity of fields.  That is, if among the generating elements of a subfield $\overline{K}$ of $K_n$, $m$ ($m\leq n$) are functionally independent, the question is reduced to proving the identity (more specifically: isomorphism) of $\overline{K}$ and $K_m$.
        \item \textit{Simplest solution of equations of multiple variables.} Let $K$ be a field of rational functions in $n$ variables $x_1,x_2,...,x_n$, which satisfy the algebraic relation $f(x_1,x_2,...,x_n) = 0$ (the case of several relations can very easily be reduced to this case).  One wishes to choose new generating elements $y_1,y_2,..,y_n$ for $K$ in such a way that the equation relating the $y_i$ is of as small of a degree as possible in one of them, say $y_1$.  What can one say about this degree?\\
        This question was treated comprehensively by Enriques \cite{ref19} at the first international congress of mathematicians (Z{\"u}rich, 1897).  I will allow myself to recall some of the nice results considered there.
        \\
        1) We consider $x_1,x_2$ as variables and $x_3,x_4,...,x_n$ as parameters. If the genus of the equation $f(x_1,x_2) = f(x_1,x_2,x_3,...,x_n) = 0$ is identically zero, one can choose a variable $t$ in such a way that $x_1,x_2$ may be expressed rationally by $t,x_3,x_4,...,x_n$ and a square root of a rational function of $x_3,x_4,...,x_n$ (M. Noether, \cite{ref56}).
        \\
        For $n = 3$, one can eliminate the quadratic irrationality using an appropriately chosen $t$.\\
        One can suspect that this theorem is meaningful for finding fields with prescribed Galois group $\mathfrak{G}$.  The latter problem seems easier to solve when the group $\mathfrak{K}/\mathfrak{G}$ is of odd order, where $\mathfrak{K}$ denotes the \textit{holomorph} of the group $\mathfrak{G}$.
        \\
        2) If the genus $p$ of the equation $f(x_1,x_2) = 0$ is $> 1$, the field $K$ can be solved by adjoining an irrationality of degree $\leq 2p -2$.\\
        Solely for the case $p = 1$ it is impossible to specify the upper bound for the degree of this irrationality.
        \item \textit{``True transcendence degree'' of an extension field}.  Let $k$ be a field of rational functions in $u_1,u_2,...,u_n$, that may be linked by an algebraic relation
        \begin{equation}\label{5.2}
            f(u_1,u_2,...,u_n) = 0.
        \end{equation}
        In addition, let $K$ be an extension field of the same transcendence degree, whose generating elements $x_1,x_2,...,x_n$ are related to the $u_i$ by the equations
        \begin{equation}\label{5.3}
            \varphi_1(x_1,x_2,...,x_n) = u_1, \quad \varphi_2(x_1,x_2,...,x_n) = u_2, ..., \varphi_n(x_1,x_2,...,x_n) = u_n.
        \end{equation}
        One should find a subfield $K_1$ of $K$ such that
        \begin{enumerate}[label = \roman*)]
            \item The compositum of $K_1$ and $k$ is precisely the field $K$;
            \item THe transcendence degree of $K_1$ is as small as possible.
        \end{enumerate}
        We will refer to the transcendence degree of $K_1$ as the true transcendence degree of $K/k$.  It is evident that every problem in field theory becomes significantly easier upon reducing the transcendence degree of the field in question.  Therein lies the significance of this problem.
        \\
        A special case of this problem manifests itself as the so-called resolvent problem (\S4).  To turn this into the resolvent problem, we must first specialize the equations (\ref{5.3}) in the following way:
        \[x_1 + x_2 + ... + x_n = -u_1, \quad x_1x_2 + ... x_{n-1}x_n = u_2, \quad x_1x_2...x_n = (-1)^nu_n.\]
        In other words, the $x_i$ are the roots of the equation
        \[x^n + u_1x^{n-1} + ... + u_{n-1}x + u_n = 0,\]
        while $k$ is generated by $u_1,u_2,...,u_n$ and a function $\Phi$ belonging to a given permutation group $\mathfrak{G}$, which itself is connected to the $u_i$ through an easily presentable equation.\\
        The resolvent problem and its extension, the problem concerning the chains of resolvents are both tasks at the limits of classical Galois theory.  The resolvent problem admits a solution via the theory of continuous groups (\S4).
        \item \textit{Rationality questions for the periods of elliptic and abelian integrals.}  We first of all consider the case of an elliptic variety. Let $K(x,y)$ be an algebraic function field, where both $x$ and $y$ are subject to the equation
        \begin{equation}\label{5.4}
            y^2 - (1-x^2)(1-k^2x^2) = 0.
        \end{equation}
        In this case it is common knowledge that one can unformize $x$ and $y$ through elliptic functions
        \begin{equation}\label{5.5}
            x = sn(u,k), \quad y = cn(u,k)dn(u,k).
        \end{equation}
        In addition let $(x_0,y_0)$ be a given point on the variety (\ref{5.4}), to which the value $u = u_0$ of the argument $u$ corresponds:
        \[x_0 = sn(u_0,k), \quad y_0 = cn(u_0,k)dn(u_0,k).\]
        $u_0$ is unambiguously determined up to a constant multiple of the periods $4K, 4K'$.  It must be decided whether $u_0$ is a rational multiple of $4K, 4K'$.
        \\
        This task may also be formulated as a problem in classical Galois theory.  As a matter of fact, we can assume that $k,x_0,y_0$ are algebraic numbers, since this case alone provides challenges.  We additionally assume that the domain of rationality $R$ contains a the module $k$.  If $u_0 \equiv \frac{mK + m'K'}{n} \text{ (modd $4K,4K'$)}$, where $m,m',n$ are integers, then $x_0 = sn\left(\frac{mK + m'K'}{n},k\right)$ satisfies a so-called divison equation
        \begin{equation}\label{5.6}
            \Phi_n(x_0,k) = 0,
        \end{equation}
        which is irreducible in the variable $k$, but can factor within $R$ if $k$ takes on certain values.  Now let $f(x_0) = 0$ be the irreducible equation satisfied by $x_0$.  The question is whether there are values of $n$ for which $\Phi_n(x,k)$ is divisible by $f(x)$.
        \item If we designate points of the Riemann surface, for which $sn\;u = sn\;u_0$ holds as $P_1$ and $P_2$, our condition can be viewed in the following sense:
        \[n\;u(P_1) -n\;u(P_2) \equiv - \text{ (modd $4K,4K'$)},\]
        where $u(P)$ is the integral of the first kind.  It follows from Abel's theorem that there exists a function $\varphi(x,y)$ belonging to the field $K(x,y)$, that has a single zero of order $n$ at $P_1$ and a single pole of order $n$ at $P_2$.  One can also say that the ideal $\frac{P_1^n}{P_2^n}$ is principal.  It follows that one can express $\varphi$ using a prime function, i.e. in the form $\varphi = e^{n\cdot\Pi(P_1,P_2)}$, where $\Pi(P_1,P_2)$ denotes the integral of the third kind with residues $+2\pi i$ at $P_1$ and $-2\pi i$ at $P_2$.
        \\
        One can construct the function $z$ of second order, that becomes infinite at $P_1$ and $P_2$ \bigg(for example, one sets due to Zolotarev $z = \frac{sn^2(u,k)}{sn^2(u,k) = sn^2(u_0,k)}$\bigg).  This brings $\Pi(P_1,P_2)$ to the form
        \[\int \frac{z+A}{\sqrt{z(z-1)(z-\alpha)(z-\beta)}}dz,\]
        where $\alpha = \frac{1}{dn^2u_0}$ and $\beta = \frac{1}{cn^2u_0}$, and it is of interest to recognize whether this integral (for suitable values of $A$, whose choice comes from the normalizing of the periods in $\Pi(P_1,P_2)$) is integrable by logarithms.
        \item This question was initially posed by Abel \cite{ref1}.  Abel solved this problem in an algebraic way, by using the fact that $\varphi = \frac{P_1^n}{P_2^n}$ is a \textit{functional unit}, i.e. the norm of $\varphi$ is constant, and solved the Diophantine equation $p^2 - q^2R$ that arises from $\varphi = \frac{p + q\sqrt{R}}{p-q\sqrt{R}}$ using the continued fraction expansion of $\sqrt{R(x)}$.  For this to be possible it is necessary and sufficient that the continued fraction expansion of $\sqrt{R(x)}$ is periodic.  However, this periodicity cannot be confirmed or denied after a finite number of steps without further effort.  To do this is was necessary to give an upper bound on the possible period.  This task was accomplished by P. Chebyshev \cite{ref77} for the case when the coefficients of $R(x)$ are rational and by G. Zolotarev \cite{ref92} for real $R(x)$ in generality.
        \\
        Chebyshev and Zolotarev produced this result, by performing the following transformations on the variables:
        \begin{equation}
            u\to u + u_0, k\to k, u_0 \to 2u_0, \alpha \to \left(\frac{\beta + \alpha - 1}{1 + \beta - \alpha}\right)^2, \beta \to \left(\frac{\beta + \alpha - 1}{1 + \alpha - \beta}\right)^2,\tag{I}
        \end{equation}
        \begin{equation}
            u\to (1 + k')u, k\to \frac{1-k'}{1+k'}, \alpha \to \left(\frac{\sqrt{(\beta - \alpha)(\beta - 1)} - \sqrt{\alpha}}{1 + \alpha - \beta}\right)^2, \beta \to \left(\frac{\sqrt{(\beta - \alpha)(\beta - 1)} + \sqrt{\alpha}}{1 + \alpha - \beta}\right)^2,\tag{II}
        \end{equation}
        over and over.  The transformation (II) has the purpose of mapping an even $n$ to an odd $n$, and is nothing more than the Landen transformation.  If after several steps the quantity $\frac{2\sqrt{k}}{1 + k}$ no longer lies in the field $K(\alpha,\beta)$, this is an indication that the integral transformed in this way corresponds to an odd $n$.  At this point one must use the transformation (I).  The answer is positive if the sequence of these transformations produces a period.  On the other hand a finite number of steps produces divisibility criteria that render the task impossible.
        \item The question regarding the commensurability of elliptic integrals is applicable to a far removed branch of mathematics, namely Diophantine analysis.  Let the equation
        \begin{equation}\label{5.7}
            f(x,y) = 0
        \end{equation}
        be given, whose coefficients belong to an algebraic number field $K$ and whose genus is $p = 1$.  We are concerned with the existence of values $x_0,y_0$ contained in the field $K$ and satisfying equation (\ref{5.7}).  In short these are called rational points of the curve (\ref{5.7}).  One can transform equation (\ref{5.7}) into the form
        \begin{equation}\label{5.8}
            y^2 = 4x^3 - g_2x - g_3,
        \end{equation}
        possibly by extending the field $K$.  The curve (\ref{5.8}) can be represented parametrically in the following manner:
        \begin{equation}\label{5.9}
            x = p(u), \quad \quad y = p'(u).
        \end{equation}
        If values $u_1, u_2$ of the argument $u$ correspond to rational points of (\ref{5.8}), it follows by the addition theorem that $u_1 \pm u_2$ also correspond to rational points.  In other words, the arguments of the rational points of the curve (\ref{5.8}) form an \textit{additive module}, that we will refer to as $K$-\textit{rationality module} in the sequel.  H. Poincar{\'e} \cite{ref59} suspected and L.J Mordell \cite{ref51} proved that each rationality module possesses a finite basis.\\
        Given $(u_1,u_2,...,u_n)$ such a basis, the question regarding the structure of the additive abelian group $\mathfrak{U}$ generated by this basis arises.  Since one can reduce every value of the argument $u$ modulo $\omega_1,\omega_2$, where $\omega_1,\omega_2$ are the periods of $p(u),p'(u)$, it follows that the number of independent generating elements of this group of infinite order (\textit{rank}) is $n$ if and only there exists no congruence of the form
        \begin{equation}\label{5.10}
            m_1u_1 + m_2u_2 + ... + m_nu_n \equiv 0 \text{(mod $\omega_1,\omega_2$)}
        \end{equation}
        between the $u_i$, where the $m_i$ are integers.  However, in general the basis of $\mathfrak{U}$ consists of a number $\rho$ of elements $A_1,A_2,...,A_{\rho}$ of infinite order, and a number of elements $B_1,B_2,...,B_{\sigma}$ ($\rho + \sigma \leq n$) of finite order.  The number $\rho$ was called the \textit{rank of rationality} of the system $u_1,u_2,...,u_n \text{ (mod $\sigma_1,\sigma_2$)}$ by Kronecker, i.e. the largest number of rationally independent elements.  The Zolotarevian process allows us to determine rank of rationality in the case $n = 1$\footnote{\textit{T. Nagell} (\cite{ref52}, p. 96, line 10-11 v.o.) says, \guillemotleft on n'a pas de m{\'e}thode g{\'e}n{\'e}rale pour reconna{\^i}tre si l'argument d'un point donn{\'e} $(x,y)$ est commensurable avec une p{\'e}riode ou non.\guillemotright [English: There is no general method to know if the argument of a given point $(x,y)$ is commensurable with a period or not.] This problem has actually solved in the aforementioned investigations by Zolotarev.}, and hopefully allows for a direct extension to the general case.
        \item One can formulate the problems arising from these reflections in the following way:
        \begin{enumerate}[label = \Roman*.]
            \item One finds a method to decide whether the value $a$ of the argument $u$ defined by the equations $sn\;a = x_0, cn\;a.dn\;a = y_0$ (or: $p(a) = x_0, p'(a) = y_0)$ is commensurable with the periods or not.\\
            It is only important to answer the question in the case where the Landen transformation of a given function is periodic.  The case is not applicable when the curve (\ref{5.4}) (or (\ref{5.8})) is real.
            \item One determines the rank of rationality of a module/modulus $(u_1,u_2,...,u_n)(\text{mod } \omega_1,\omega_2)$, where the $u_i$ are to be determined from the equations $p(u_i) = x_i, p'(u_i) = y_i$ and $\omega_1,\omega_2$ denote the periods of the elliptic function $p(u)$.
            \item Let an elliptic function field (perhaps through specification of the module/modulus $k$) and an algebraic number field $K$ containing $k$ be given.  One must find a basis of the corresponding $K$-rationality module/modulus.  If a basis $(u_1,u_2,...,u_n)$ (through explicit specification of the values of $p(u_i))$ is given, one decided whether it is a basis of the $K$-rationality module/modulus or not.
        \end{enumerate}
        \item One can transfer the formulation of the aforementioned tasks to general algebraic functions without much fuss.  First of all it seems reasonable that the question regarding the finiteness of an algebraically specified transformation of the ``Jacobian group'' descibed in \S1, No. 9 is closely related to the question regarding the integrability of abelian integrals by logarithms.  On the other hand one can readily transfer the reflections in No. 10 to the general case.  This also shows that there is a certain relationship between the Jacobian group and the ideal class groups of an algebraic function field.  However, this problem is linked with the question concerning the functional units.\footnote{This connection was brought to my attention by my highly respected teacher Prof. Dr. D. Grave.  See for example Verh. Russ. Math. Kongre{\ss} in Moscow (1927), p. 215 (Russian).}\\
        A. Weil \cite{ref87} investigated Diophantine equations $f(x,y) = 0$ of arbitrary $p$ using a similar method.  In particular he extended Mordell's resullt about the finiteness of $K$-rationality modules/moduli for arbitrary $p \geq 1$.  To do so he used the Jacobian group, while Mordell \cite{ref51} implicitly used the division of elliptic functions.
        \item The Diophantine problem can be considered as a special case of the Hilbert-Doerge irreducibility problem:
        \\Let an equation $f(z,t) = 0$ be given.  One must find all values $t_i$ of $t$ for which $f(z,t_i)$ is a reducible polynomial in a prescribed number field $k$.
        \\
        According to $K$ Doerge it follows from the investigations of Weil \cite{ref87} that for $p > 1$ $t$ takes on infinitely many such values only if $f(z,t)$ can be transformed using a transformation of the form
        \[t = c_{-m}u^{-m} + c_{-m+1}u^{-m+1} +... + c_0 +c_1u + ... + c_mu^m \]
        into an identically decomposing polynomial in $z$ and $u$ \footnote{Amendment while correcting.  In this direction significant new results have been obtained by C.R. Siegel (Abh. preuss. Akad., Berlin, 1930, No 1).}.
        \\
        The work of Doerge under discussion is devoted to the case where $z$ is a functional unit.  Then Doerge obtains very easy conditions for when $f(z,t)$ is only reducible for a finite number of values of $t$ in $k$.  It is remarkable that this results in a new direct connection between the Diophantine equations and the functional units.
        \\
        \end{enumerate}
(Submitted the 12th of September, 1932)
\newpage
\bibliographystyle{acm}
\bibliography{main}

\begin{thebibliography}{10}

\bibitem{ref1}
{\sc Abel, N.}
\newblock Sur l'int{\'e}gration de la formule etc.
\newblock In {\em Oeuvres Completes, Tome {I}}. Christiania {I}mprimerie De
  Grondahl \& Son, 1881, ch.~11, pp.~104--144.

\bibitem{ref2}
{\sc Artin, E.}
\newblock Beweis des allgemeinen {R}eziprozit{\"a}tsgesetzes.
\newblock {\em Hamb. Abh. 5\/} (1927), 353--363.

\bibitem{ref3}
{\sc Baer, R.}
\newblock Abbildungseigenschaften algebraischer {E}rweiterungen.
\newblock {\em Math. Zeitschrift 33\/} (1931), 451--479.

\bibitem{ref4}
{\sc Baker, H.}
\newblock {\em Abel's {T}heorem and the allied {T}heory etc.}
\newblock Cambridge, 1897.

\bibitem{ref5}
{\sc Bauer, M.}
\newblock Ganzzahlige {G}leichungen ohne {A}ffekt.
\newblock {\em Math. Ann. 64\/} (1907), 325--327.

\bibitem{ref6}
{\sc Bauer, M.}
\newblock Ganzzahlige {G}leichungen ohne {A}ffekt.
\newblock {\em Math. Zeitschrift 16\/} (1923), 318--319.

\bibitem{ref7}
{\sc Brandt, H.}
\newblock Ueber eine {V}erallgemeinerung des {G}ruppenbergriffes.
\newblock {\em Math. Ann. 96\/} (1926), 360--366.

\bibitem{ref8}
{\sc Brauer, R.}
\newblock {\em Ueber die {D}arstellung der {D}rehungsgruppe durch {G}ruppen
  linearer {S}ubstitutionen}.
\newblock PhD thesis, University of Berlin, 1925.

\bibitem{ref9}
{\sc Brauer, R., Hasse, H., and Noether, E.}
\newblock Beweis eines {H}auptsatzes in der {T}heorie der {A}lgebren.
\newblock {\em Crelle (Hensel-Festband) 167\/} (1931), 399--404.

\bibitem{ref10}
{\sc Breuer, S.}
\newblock Zur {B}estimmung der metazyklischen {M}inimalbasis vom
  {P}rimzahlgrad.
\newblock {\em Math. Ann. 92\/} (1924).

\bibitem{ref11}
{\sc Breuer, S.}
\newblock Metazyklische {M}inimalbasis und komplexe {P}rimzahlen.
\newblock {\em Crelle 156\/} (1927), 13--42.

\bibitem{ref12}
{\sc Bring, E.~S.}
\newblock Meletamata quaedam mathematica circa transformationem etc.
\newblock {\em Diss. Lund\/} (1786).

\bibitem{ref13}
{\sc Cartan, {\'E}.}
\newblock {\em Sur la structure des groupes finis et continus}.
\newblock PhD thesis, University of Paris, {{\'E}}cole {N}ormale
  {S}up{\'e}rieure, 1894.

\bibitem{ref14}
{\sc Cartan, {\'E}.}
\newblock Sur la structure des groupes infinis.
\newblock {\em C R. 135\/} (1902), 851--854.

\bibitem{ref15}
{\sc Castelnuovo, G.}
\newblock Sulla razionalit{\`a} delle involuzioni piane.
\newblock {\em Math. Ann. 44\/} (1894), 125--155.

\bibitem{ref16}
{\sc Dedekind, R.}
\newblock Zur {T}heorie der {{I}}deale.
\newblock {\em G{\"o}tt. Nachr.\/} (1894), 272--277.

\bibitem{ref17}
{\sc Dedekind, R.}
\newblock Ueber die {A}nzahl der {{I}}dealklassen in reinen kubischen
  {Z}ahlk{\"o}rpern.
\newblock {\em Crelle 121\/} (1900), 40--123.

\bibitem{ref18}
{\sc Doerge, K.}
\newblock Bemerkung zum {H}ilbertschen {{I}}rreduzibilit{\"a}tssatz.
\newblock {\em Math. Ann. 102\/} (1929), 521--530.

\bibitem{ref19}
{\sc Enriques, F.}
\newblock Sur les probl{\`e}mes qui se rapportent {\`a} la r{\'e}solution des
  {\'e}quations alg{\'e}briques etc.
\newblock {\em Math. Ann. 51\/} (1899), 134--153.

\bibitem{ref20}
{\sc Enriques, F.}
\newblock Sopra una involuzione non razionale dello spazio.
\newblock {\em Rendic. Linc. 21\/} (1912), 81--83.

\bibitem{ref21}
{\sc Fano, G.}
\newblock Sopra alcune variet{\`a} algebriche a tre dimensioni etc.
\newblock {\em Atti Acc. Torino 43\/} (1908), 973--981.

\bibitem{ref22}
{\sc Fischer, E.}
\newblock Zur {T}heorie der endlichen {A}belschen {G}ruppen.
\newblock {\em Math. Ann. 77\/} (1916), 81--88.

\bibitem{ref23}
{\sc Frobenius, G.}
\newblock Ueber {B}eziehungen zwischen {P}rimidealen eines algebraischen
  {K}{\"o}rpers und den {S}ubstitutionen usw.
\newblock {\em Sitzber. Berl. Akad.\/} (1896), 689--705.

\bibitem{ref24}
{\sc Fueter, R.}
\newblock Die {T}heorie der {Z}ahlstrahlen {I},{I}{I}.
\newblock {\em Crelle 130,132\/} (1905,1907), 197--257,255--269.

\bibitem{ref25}
{\sc Fueter, R.}
\newblock Abelsche {G}leichungen in quadratisch imagin{\"a}ren
  {Z}ahlk{\"o}rpera.
\newblock {\em Math. Ann. 75\/} (1914), 177--255.

\bibitem{ref26}
{\sc Furtw{\"a}ngler, P.}
\newblock Allgemeiner {E}xistenzbeweis f{\"u}r den {K}lassenk{\"o}rper eines
  beliebigen {Z}ahlk{\"o}rpers.
\newblock {\em Math. Ann. 63\/} (1907), 1--37.

\bibitem{ref27}
{\sc Furtw{\"a}ngler, P.}
\newblock Ueber {M}inimalbasen f{\"u}r {K}{\"o}rper rationaler {F}unktionen.
\newblock {\em Sitzber. Wiener Akad. 134\/} (1925), 69--80.

\bibitem{ref28}
{\sc Garver, R.}
\newblock On the removal of four terms from an equation by means of a
  {T}schirnhaus transformation.
\newblock {\em Bull. Amer. Math. Soc. 35\/} (1929), 73--78.

\bibitem{ref29}
{\sc Halphen, G.~H.}
\newblock Trait{\'e} des fonctions {\'e}lliptiques et leurs applications etc.
\newblock {\em Tome 3. Paris\/} (1891).

\bibitem{ref31}
{\sc Hasse, H.}
\newblock Ein weiteres {E}xistenztheorem in der {T}heorie der algebraischen
  {Z}ahlk{\"o}rper.
\newblock {\em Math. Zeitschrift 24\/} (1925), 149--160.

\bibitem{ref30}
{\sc Hasse, H.}
\newblock Zwei {E}xistenztheoreme {\"u}ber algebraische {Z}ahlk{\"o}rper.
\newblock {\em Math. Ann. 95\/} (1925), 229--238.

\bibitem{ref32}
{\sc Hasse, H.}
\newblock Bericht {\"u}ber neuere {U}ntersuchungen und {P}robleme aus der
  {T}heorie der algebraischen {Z}ahlk{\"o}rper.
\newblock {\em {I}, Jahresber. D. M. V.; {I}a, ibid. 26; {I}{I}, ibid.; V{I}
  {E}rg{\"a}nzbd. 35\/} (1926,1927,1930), 1--55, 233--311.

\bibitem{ref33}
{\sc Hasse, H.}
\newblock Arithmetische {T}heorie der kubischen {Z}ahlk{\"o}rper auf
  klassenk{\"o}rpertheoretischer {G}rundlage.
\newblock {\em Math. Zeitschrift 31\/} (1930), 565--582.

\bibitem{ref34}
{\sc Hasse, H.}
\newblock Theory of {C}yclic {A}lgebras over an {A}lgebraic {N}umber {F}ield.
\newblock {\em Trans. Amer. Math. Soc. 34\/} (1932), 171--214.

\bibitem{ref35}
{\sc Hilbert, D.}
\newblock Ueber die irreducibilit{\"a}t ganzer rationaler {F}unktionen usw.
\newblock {\em Crelle 110\/} (1892), 104--129.

\bibitem{ref36}
{\sc Hilbert, D.}
\newblock Mathematische {P}robleme.
\newblock {\em G{\"o}tt. Nachr.\/} (1900), S. 253--297.

\bibitem{ref37}
{\sc Hilbert, D.}
\newblock Ueber die {G}leichung neunten {G}rades.
\newblock {\em Math. Ann. 97\/} (1926), 243--250.

\bibitem{ref38}
{\sc Hurwitz, A.}
\newblock Ueber algebraische {G}ebilde mit eindeutigen {T}ransformationen in
  sich.
\newblock {\em Math. Ann. 41\/} (1893), 403--442.

\bibitem{ref39}
{\sc Jordan, C.}
\newblock M{\'e}moire sur les {\'e}quations diff{\'e}rentielles etc.
\newblock {\em Crelle 84\/} (1878), 89--215.

\bibitem{ref40}
{\sc Killing, W.}
\newblock Die {Z}usammensetzung der stetigen endlichen
  {T}ransformationsgruppen. {I}, {I}{I}, {I}{I}{I}, {I}v.
\newblock {\em Math. Ann. 31 33, 34, 36\/} (1888), 252--290, 1--48, 57--122,
  161--189.

\bibitem{ref41}
{\sc Klein, F.}
\newblock Gesammelte mathematische {A}bhandlungen, berlin, bd. 2.
\newblock 255--504.

\bibitem{ref42}
{\sc Kronecker, L.}
\newblock Vorlesungen {\"u}ber {Z}ahlentheorie.
\newblock {\em Lpz.\/} (1901), 452--492.

\bibitem{ref43}
{\sc Krull, W.}
\newblock Galoissche {T}heorie der unendlichen algebraischen {E}rweiterungen.
\newblock {\em Math. Ann. 100\/} (1928), 687--698.

\bibitem{ref44}
{\sc Levi, E.~E.}
\newblock Sulla struttura dei gruppi finiti e continui.
\newblock {\em Atti Acc. Torino 40\/} (1905), 423--437.

\bibitem{ref45}
{\sc Lie, S.}
\newblock Theorie der {T}ransformationsgruppen, bd. 1.
\newblock {\em Lpz.\/} (1888).

\bibitem{ref46}
{\sc Loewy, A.}
\newblock Neue elementare {B}egr{\"u}ndung und {E}rweiterung der {G}aloisschen
  {T}heorie.
\newblock {\em Sitzber. Heidlb. Akad. {I}, 7. Abh 1925; {I}{I}, {I} Abh.
  1927\/}.

\bibitem{ref47}
{\sc Loewy, A.}
\newblock Ueber abstrakt definierte {T}ransmutationssysteme oder
  {M}ischgruppen.
\newblock {\em Crelle 157\/} (1927), 239--254.

\bibitem{ref48}
{\sc {L{\"u}roth}}.
\newblock Beweis eines {S}atzes {\"u}ber rationale {C}urven.
\newblock {\em Math. Ann. 9\/} (1876), 163--165.

\bibitem{ref49}
{\sc Mertens, F.}
\newblock Ueber {D}irichlet's {B}eweis des {S}atzes, da{\ss} jede unbegrenzte
  arithmetische {P}rogression, deren {D}ifferenz zu ihren {G}liedern
  teilerfremd ist, unendlich viele {P}rimzahlen enth{\"a}lt.
\newblock {\em Sitzber. Wiener Akad. 106\/} (1897), 254--286.

\bibitem{ref50}
{\sc Mertens, F.}
\newblock Ein {B}eweis des {G}alois'schen {F}undamentalsatzes.
\newblock {\em Sitzber. Wiener Akad. 111\/} (1902), 17--37.

\bibitem{ref51}
{\sc Mordell, L.~J.}
\newblock On the {R}ational {S}olutions of the {{I}}ndeterminate equations of
  the {T}hird and {F}ourth {D}egrees.
\newblock {\em Proc. Cambr. Phil. Soc. 21\/} (1922), 179--192.

\bibitem{ref52}
{\sc Nagell, T.}
\newblock Sur les propri{\'e}t{\'e}s arithm{\'e}tiques des cubiques planes du
  premier genre.
\newblock {\em Acta Math. 52\/} (1928), 93--126.

\bibitem{ref53}
{\sc Netto, E.}
\newblock Ueber einen {L}{\"u}roth-{G}ordan'schen {S}atz.
\newblock {\em Math. Ann. 46\/} (1895), 310--318.

\bibitem{ref54}
{\sc Neumann, C.}
\newblock {\em Vorlesungen {\"u}ber {R}iemann's {T}heorie der {A}bel'schen
  {{I}}ntegrale,}.
\newblock 2 Aufl., Lpz, 1884.

\bibitem{ref55}
{\sc Noether, E.}
\newblock Gleichungen mit vorgeschriebener {G}ruppe.
\newblock {\em Math. Ann. 78\/} (1918), 221--227.

\bibitem{ref56}
{\sc Noether, M.}
\newblock Ueber {F}l{\"a}chen, welche {S}chaaren rationaler {C}urven besitzen.
\newblock {\em Math. Ann. 3\/} (1871), 161--227.

\bibitem{ref57}
{\sc Ore, O.}
\newblock Zur {T}heorie der {E}isensteinschen {G}leichungen.
\newblock {\em Math. Zeitschrift 20\/} (1924), 267--279.

\bibitem{ref58}
{\sc Perron, O.}
\newblock Ueber {G}leichungen ohne {A}ffekt.
\newblock {\em Sitzber. Heidlb. Akad. 3 Abh.\/} (1923).

\bibitem{ref59}
{\sc Poincar{\'e}, H.}
\newblock Sur les propri{\'e}t{\'e}s arithm{\'e}tiques des courbes
  alg{\'e}briques.
\newblock {\em Journ. de Math. 17\/} (1901), 161.

\bibitem{ref60}
{\sc Pollaczek, F.}
\newblock Ueber die {E}inheiten relativ-{A}belscher {Z}ahlk{\"o}rper.
\newblock {\em Math. Zeitschrift 30\/} (1929), 520--551.

\bibitem{ref61}
{\sc Remak, R.}
\newblock Ueber die {A}bsch{\"a}tzung des absoluten {B}etrages des {R}egulators
  eines algebraischen {Z}ahlk{\"o}rpers nach unten.
\newblock {\em Crelle 167\/} (1931), 360--378.

\bibitem{ref62}
{\sc Schatunowski, S.}
\newblock {\em Algebra als {L}ehre von den {K}ongruenzen nach funktionalen
  {M}oduln (russisch)}.
\newblock PhD thesis, Odessa, 1917.

\bibitem{ref64}
{\sc Scholz, A.}
\newblock Reduktion der {K}onstruktion von k{\"o}rpern mit zweistufiger
  (metaabelscher) {G}ruppe.
\newblock {\em Sitzber. Heidlb. Akad. 14. Abh.\/} (1929).

\bibitem{ref63}
{\sc Scholz, A.}
\newblock Ueber die {B}ildung algebraischer {Z}ahlk{\"o}rper mit
  aufl{\"o}sbarer {G}aloisscher {G}ruppe.
\newblock {\em Math. Zeitschrift 30\/} (1929), 332--356.

\bibitem{ref65}
{\sc Scholz, A.}
\newblock Ein {B}eitrag zur {T}heorie der {Z}usammensetzung endlicher
  {G}ruppen.
\newblock {\em Math. Zeitschrift 32\/} (1930), 187--189.

\bibitem{ref66}
{\sc Scholz, A.}
\newblock Ueber das {V}erh{\"a}ltnis von {{I}}dealklassen-un {E}inheitengruppe
  in {A}belschen {K}{\"o}rpern von {P}rimzahlpotenzgrad.
\newblock {\em Sitzber. heidlb. Akad. 3. Abh.\/} (1930), 31--55.

\bibitem{ref67}
{\sc Schreier, O.}
\newblock Die {V}erwandtschaft stetiger {G}ruppen im gro{\ss}en.
\newblock {\em Hamb. Abh. 5\/} (1927), 233--244.

\bibitem{ref68}
{\sc Schur, I.}
\newblock {\em Ueber eine {K}lasse von {M}atrizen, die sich einer gegebenen
  {M}atrix zuordnen lassen}.
\newblock PhD thesis, Berlin, 1901.

\bibitem{ref69}
{\sc Schur, I.}
\newblock Ueber die {D}arstellung der endlichen {G}ruppen durch gebrochene
  lineare {S}ubstitutionen {I}\&{I}{I}.
\newblock {\em Crelle 127,132\/} (1904, 1907), 20--50,85--137.

\bibitem{ref70}
{\sc Schur, I.}
\newblock Beispiele f{\"u}r {G}leichungen ohne {A}ffekt.
\newblock {\em Jahresber. DMV. 29\/} (1920).

\bibitem{ref71}
{\sc Schur, I.}
\newblock Ueber die stetigen {D}arstellungen der allgemeinen linearen {G}ruppe.
\newblock {\em Sitzber. Berl. Akad.\/} (1928), 96--124.

\bibitem{ref72}
{\sc Schur, I.}
\newblock Gleichungen ohne {A}ffekt.
\newblock {\em Sitzber. Berl. Akad.\/} (1930), 443--449.

\bibitem{ref73}
{\sc Speiser, A.}
\newblock Die {Z}erlegung von {P}rimzahlen in algebraischen {Z}ahlk{\"o}rpern.
\newblock {\em Trans. Amer. Math. Soc. 23\/} (1922), 173--178.

\bibitem{ref74}
{\sc Steinitz, E.}
\newblock Algebraische {T}heorie der {K}{\"o}rper.
\newblock {\em Crelle 137\/} (1910), 167--308.

\bibitem{ref75}
{\sc Sylvester, J.~J.}
\newblock On the so-called {T}schirnhausen {T}ransformation.
\newblock {\em Crelle 100\/} (1886), 465--487.

\bibitem{ref76}
{\sc Takagi, T.}
\newblock Ueber eine {T}heorie des relativ {A}bel'schen {Z}ahlk{\"o}rpers.
\newblock {\em Journ. Coll. Sc. Tokyo 41, Art. 9\/} (1920).

\bibitem{ref77}
{\sc Tchebycheff, P.}
\newblock Sur l'int{\'e}gration de la diff{\'e}rentielle$\frac{x + A}{\sqrt{x^4
  + \alpha x^3 + \beta x^2 + \gamma x + \delta }}dx.$.
\newblock {\em Oeuvres 1\/}, S. 517--530.

\bibitem{ref78}
{\sc Tschebotar{\"o}w, N.}
\newblock Die der {T}schirnhausenschen umgekehrte {A}ufgabe (russisch).
\newblock {\em Journal des Sciences 1\/} (1922).

\bibitem{ref80}
{\sc Tschebotar{\"o}w, N.}
\newblock Die {B}estimmung der {D}ichtigkeit einer {M}enge von {P}rimzahlen,
  welche zu einer gegebenen {S}ubstitutionsklasse geh{\"o}ren.
\newblock {\em Math. Ann 95\/} (1925), 191--228.

\bibitem{ref79}
{\sc Tschebotar{\"o}w, N.}
\newblock Zur {A}ufgabe der {B}estimmung von algebraischen {G}leichungen mit
  vorgeschriebener {G}ruppe.
\newblock {\em Bull. Soc. Math. Kasan 1\/} (1926), 26--32.

\bibitem{ref81}
{\sc Tschebotar{\"o}w, N.}
\newblock Studien {\"u}ber {P}rimzahlendichtigkeiten {I} \& {I}{I}.
\newblock {\em Bull. Soc. Math. Kasan 2,3\/} (1927,1928), 1--17.

\bibitem{ref82}
{\sc Tschebotar{\"o}w, N.}
\newblock Zur {G}ruppentheorie des {K}lassenk{\"o}rpers.
\newblock {\em Crelle 161\/} (1929), 179--193.

\bibitem{ref84}
{\sc Tschebotar{\"o}w, N.}
\newblock Untersuchungen {\"u}ber relativ {A}belsche {Z}ahlk{\"o}rper.
\newblock {\em Crelle 167\/} (1931), 98--121.

\bibitem{ref83}
{\sc Tschebotar{\"o}w, N.}
\newblock Ueber ein algebraisches {P}roblem von {H}errn {H}ilbert {I} \&
  {I}{I}.
\newblock {\em Math. Ann. 104, 105\/} (1931, 1931), 459--471,240--255.

\bibitem{ref85}
{\sc Voronoi, G.}
\newblock Ueber ganze algebraische {Z}ahlen, die von einer {W}urzel einer
  {G}leichung 3. {G}rades abh{\"a}ngen (russisch).
\newblock Master's thesis, S-Pb, 1894.

\bibitem{ref86}
{\sc Weber, H.}
\newblock {\em Lehrbuch der {A}lgebra, Bd. 3}.
\newblock Braunschweig, 1908.

\bibitem{ref87}
{\sc Weil, A.}
\newblock L'arithm{\'e}tique sur les courbes alg{\'e}briques.
\newblock {\em Acta Math. 52\/} (1929), 281--315.

\bibitem{ref88}
{\sc Weyl, H.}
\newblock Theorie der {D}arstellung kontinuierlicher halb-einfacher {G}ruppen
  durch lineare {T}ransformationen {I} \& {I}{I} \& {I}{I}{I}.
\newblock {\em Math. Zeitschrift 23, 24,24\/} (1925, 1925, 1925), 271--309,
  328--376, 377--395.

\bibitem{ref89}
{\sc Wiman, A.}
\newblock Ueber eine einfache {G}ruppe von 360 ebenen {C}ollineationen.
\newblock {\em Math. Ann. 47\/} (1896), 531--556.

\bibitem{ref90}
{\sc Wiman, A.}
\newblock Ueber die {D}arstellung der symmetrischen und alternierenden
  {V}ertauschungsgruppen usw.
\newblock {\em Math. Ann. 52\/} (1899), 243--270.

\bibitem{ref91}
{\sc Wiman, A.}
\newblock Ueber die {A}nwendung der {T}schirnhausentransformation auf die
  {R}eduktion algebraischer {G}leichungen.
\newblock {\em Noa Acta Uppsala\/} (1927), X + 3--8.

\bibitem{ref92}
{\sc Zolotareff, G.}
\newblock Th{\'e}orie des nombres complexes entiers avec une application vers
  le calcul int{\'e}gral.
\newblock {\em Diss. S.-Pb. 1874. Oeuvres tome {I}, Leningrad\/} (1931),
  161--360.

\end{thebibliography}
\end{document}